\documentclass[twoside,12pt]{article}
\usepackage{geometry}
\usepackage{amsmath,amssymb,amsfonts,amsthm}
\usepackage{txfonts}
\usepackage{todonotes}
\usepackage{tablists}
\usepackage{graphicx}
\usepackage{mathrsfs}
\usepackage{color}
\usepackage{float}
\usepackage{extarrows}
\usepackage{exscale}
\usepackage{titlesec}
\titleformat{\section}[block]{\large\center\sc}{\arabic{section}}{0.5em}{}[] 
\usepackage{relsize}
\usepackage{color}
\usepackage{hyperref}
\usepackage[titles]{tocloft} 
\usepackage{fancyhdr}
\theoremstyle{plain}
\newtheorem{theorem}{Theorem}[section]
\newtheorem{lemma}[theorem]{Lemma}

\newtheorem{proposition}[theorem]{Proposition}

\newtheorem{definition}[theorem]{Definition}

\newtheorem{problem}[theorem]{Problem}

\newtheorem{remark}[theorem]{Remark}

\let\oldsection\section
\renewcommand\section{\setcounter{equation}{0}\oldsection}

\def\be{\begin{equation}}
\def\ee{\end{equation}}
\def\bes{\begin{equation*}}
\def\ees{\end{equation*}}
\def\bs{\begin{split}}
\def\es{\end{split}}
\def\bali{\begin{aligned}}
\def\eali{\end{aligned}}
\newcommand{\pf}{\noindent {\bf Proof. \hspace{2mm}}}
\def\bR{{\mathbb R}}
\def\un{\underbrace}
\def\al{\alpha}

\def\ve{\varepsilon}

\def\g{\gamma}

\def\Dl{\Delta}
\def\lt{\left}
\def\rt{\right}
\def\ls{\lesssim}

\def\i{\infty}
\def\p{\partial}
\def\f{\frac}
\def\na{\nabla}

\def\O{\Omega}

\def\q{\quad}

\def\bl{\boldsymbol}
\def\mS{\mathbb{S}}
\def\mR{\mathbb{R}}

\def\mH{\mathcal{H}}

\allowdisplaybreaks
\def\mD{\mathcal{D}}
\def\mfD{\mathfrak{D}}

\begin{document}

\title{\bf\normalsize ON LERAY'S PROBLEM IN AN INFINITE LONG PIPE WITH THE NAVIER-SLIP BOUNDARY CONDITION}

\author{\normalsize\sc Zijin Li, Xinghong Pan and Jiaqi Yang}

\date{}

\maketitle

\begin{abstract}
The original Leray's problem concerns the well-posedness of weak solutions to the steady incompressible Navier-Stokes equations in a distorted pipe, which approach to the Poiseuille flow subject to the no-slip boundary condition at spacial infinity. In this paper, the same problem with the Navier-slip boundary condition instead of the no-slip boundary condition, is addressed. Due to the complexity of the boundary condition, some new ideas, presented as follows, are introduced to handle the extra difficulties caused by boundary terms.

First, the Poiseuille flow in the semi-infinite straight pipe with the Navier-slip boundary condition will be introduced, which will be served as the asymptotic profile of the solution to the generalized Leray's problem at spacial infinity. Second, a solenoidal vector function defined in the whole pipe, satisfying the Navier-slip boundary condition, having the designated flux and equalling to the Poiseuille flow at large distance, will be carefully constructed. This plays an important role in reformulating our problem. Third, the energy estimates depend on a combined $L^2$ estimate of the gradient and the stress tensor of the velocity.

\medskip



\end{abstract}

\tableofcontents

\section{Introduction}
\q\ The 3D stationary Navier-Stokes (NS) equations which describes the motion of stationary viscous incompressible fluids follows that
\be\label{NS}
\lt\{
\begin{aligned}
&\bl{u}\cdot\na \bl{u}+\na p-\Dl \bl{u}=0,\\
&\na\cdot \bl{u}=0,
\end{aligned}
\rt.\q \text{in}\q \mathcal{D}\subset \bR^3.
\ee
Here  $\bl{u}(x)\in\mathbb{R}^3$, $p(x)\in\mathbb{R}$ represent the velocity and the scalar pressure respectively. In this paper, we consider the smooth domain $\mD$ to be an infinitely long pipe with two straight ``outlets", the left and the right, and a compact distortion, a bubble or a bulge, in the middle. We denote it as follows.

\be\label{cylinder}
\mD=\mD_L\cup\mD_M\cup\mD_R,
\ee
where $\mD_L$ and $\mD_R$ are semi-infinite smooth straight pipes with their cross sections $\Sigma_L$ and $\Sigma_R$ being perpendicular to $x_3$-axis, namely:
\[
\mD_L=\Sigma_L\times(-\i, -Z/2],\q\text{and}\q\mD_R=\Sigma_R\times[Z/2,+\i).
\]
Here $\Sigma_L,\,\Sigma_R\subset\mR^2$ are smooth bounded domains. The distortion part $\mD_M\subset\mR^2\times(-Z,Z)$ is a compact smooth domain in $\bR^3$.

Moreover, technically we assume there exists an infinitely long smooth straight pipe getting through $\mD$, which means there exists $\Sigma'\subset\subset\Sigma_L\cap\Sigma_R$ such that $\Sigma'\times\mR\subset\mD$. This will be applied in the construction of the profile vector in Section \ref{SEC3.2}.
\begin{figure}[H]
\centering
\includegraphics[scale=0.5]{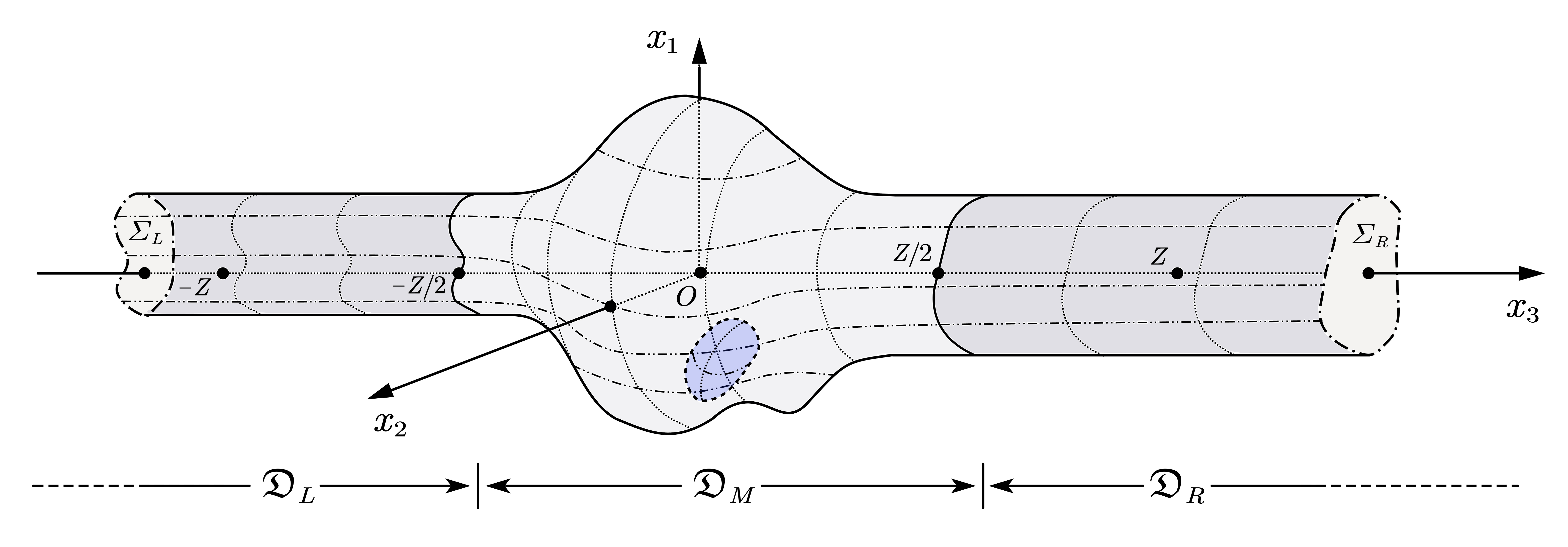}
\caption{Infinite pipe $\mD$, with a bubble and an obstacle in the middle}
\end{figure}
In the current paper, the Navier-Stokes equations \eqref{NS} will be equipped with the following boundary condition:
{\noindent\bf The Navier-slip boundary condition}:
\be\label{NBC}
\left\{
\begin{aligned}
&2(\mathbb{S}\bl{u}\cdot\boldsymbol{n})_{\mathrm{tan}}+\alpha \bl{u}_{\mathrm{tan}}=0,\\
&\bl{u}\cdot\boldsymbol{n}=0,\\
\end{aligned}
\right.\q\text{on}\q\p\mD.
\ee
Here $\mathbb{S}\bl{u}=\frac{1}{2}\left(\nabla \bl{u}+\nabla^T \bl{u}\right)$ is the stress tensor, where $\nabla^T \bl{u}$ stands for the transpose of the Jacobian matrix $\nabla \bl{u}$, and $\boldsymbol{n}$ is the unit outer normal vector of $\p \mD$. For a vector field $\bl{v}$, we denote $\bl{v}_{\mathrm{tan}}$ its tangential part:
\[
\bl{v}_{\mathrm{tan}}:=\bl{v}-(\bl{v}\cdot\boldsymbol{n})\boldsymbol{n}.
\]
And $\al>0$ stands for the friction constant which may depend on various elements, such as the property of the boundary and the viscosity of the fluid. When $\al\to0_+$, the boundary condition \eqref{NBC} turns to be the total Navier-slip boundary condition, while when $\al\to\infty$, the boundary condition \eqref{NBC} degenerates into the no-slip boundary condition $u\equiv 0$ on the boundary. In this paper, we will assume $0<\al<+\i$.



Throughout this paper, $C_{a,b,c,...}$ denotes a positive constant depending on $a,\,b,\, c,\,...$, which may be different from line to line. For a vector $x=(x_1,x_2,x_3)\in\mR^3$, we denote $x_h:=(x_1,x_2)$. For a two-dimensional scalar function $f$ or a vector-valued function $\bl{f}:=(f_1,f_2)$, we denote
\[
\nabla_h f=(\p_{x_1}f,\p_{x_2}f),\q \Dl_hf=\p^2_{x_1}f+\p^2_{x_2}f,\q\text{div}_h\bl{f}:=\p_{x_1}f_1+\p_{x_2}f_2.
\]
Meanwhile, for any $\zeta>1$, we denote
\[
\mD_\zeta:=\left\{x\in\mathcal{D}:\,\,-\zeta<x_3<\zeta\right\},
\]
the truncated pipe with the length of $2\zeta$. Meanwhile, notations $\O^{\pm}_\zeta$ is denoted by
\[
\O^{+}_\zeta:=\left(\mathcal{D}_\zeta-\mathcal{D}_{\zeta-1}\right)\cap\{x\in\mathcal{D}:\,x_3>0\},\q\O^{-}_\zeta:=\left(\mathcal{D}_\zeta-\mathcal{D}_{\zeta-1}\right)\cap\{x\in\mathcal{D}:\,x_3<0\},
\]
respectively. We also apply $A\lesssim B$ to state $A\leq CB$. Moreover, $A\simeq B$ means both $A\lesssim B$ and $B\lesssim A$.

For $1\leq p\leq\infty$ and $k\in\mathbb{N}$, $L^p$ denotes the usual Lebesgue space with norm
\[
\|f\|_{L^p(D)}:=
\lt\{
\begin{aligned}
&\left(\int_{D}|f(x)|^pdx\right)^{1/p},\quad &1\leq p<\infty,\\[3mm]
&\mathrm{ess sup}_{x\in D}|f(x)|,\quad &p=\infty,\\
\end{aligned}
\rt.
\]
while $W^{k,p}$ denotes the usual Sobolev space with its norm
\[
\begin{split}
\|f\|_{W^{k,p}(D)}:=&\sum_{0\leq|L|\leq k}\|\nabla^L f\|_{L^p(D)},\\
\end{split}
\]
where $L=(l_1,l_2,l_3)$ is a multi-index. We also simply denote $W^{k,p}$ by $H^k$ provided $p=2$. Finally, $\bar{D}$ denote the closure of a domain $D$. A function $g\in W^{k,p}_{\mathrm{loc}}(D)$ or $W^{k,p}_{\mathrm{loc}}(\bar{D})$ function means $g\in W^{k,p}(\tilde{D})$, for any $\tilde{D}$ compactly contained in $D$ or $\bar{D}$.

For the 3D vector-valued function, we define
\[
\begin{split}
\mathcal{H}(\mD)&:=\left\{\bl{\varphi}\in H^1({\mD};\bR^3):\,\bl{\varphi}\cdot\bl{n}\big|_{\p \mD}=0 \right\},\\
\mathcal{H}_\sigma({\mD})&:=\left\{\bl{\varphi}\in H^1({\mD};\,\mathbb{R}^3):\,\na\cdot \bl{\varphi}=0,\ \bl{\varphi}\cdot \bl{n}\big|_{\p \mD}=0\right\},
\end{split}
\]
and
\[
\mathcal{H}_{\sigma,{\mathrm{loc}}}(\overline{\mD}):=\left\{\bl{\varphi}\in H^1_{{\mathrm{loc}}}(\overline{\mD};\,\mathbb{R}^3):\,\na\cdot \bl{\varphi}=0,\ \bl{\varphi}\cdot \bl{n}\big|_{\p \mD}=0\right\}.
\]
We also denote
\[
 \bl{X}:=\big\{\bl{\varphi}\in C_c^\infty(\overline{\mD}\,;\,\mathbb{R}^3):\ \na\cdot \bl{\varphi}=0,\ \bl{\varphi}\cdot\bl{n}\big|_{\p\mD}=0\big\}.
\]
Clearly, $\bl{X}$ is dense in $\mathcal{H}_\sigma$ in $H^1(\mD)$ norm. For matrices $\bl{\Gamma}=(\gamma_{ij})_{1\leq i,j\leq 3}$ and $\bl{K}=(\kappa_{ij})_{1\leq i,j\leq 3}$, we denote
\[
\bl{\Gamma}:\bl{K}=\sum^{3}_{i,j=1}\gamma_{ij}\kappa_{ij}.
\]
Next, we state our main problem of the paper.
\subsection{The Leray's problem with the Navier-slip boundary condition}
\q\ For a given flux $\Phi$ which is supposed to be nonnegative without loss of generality, if we consider the Poiseuille  flow, $\boldsymbol{g}^i_{\Phi}$, of \eqref{NS} with the boundary condition \eqref{NBC} in $\mD_i$ ($i$ denotes $L$ or $R$), then it satisfies
\[
\left\{
\begin{aligned}
&\boldsymbol{g}^i_{\Phi}=g^i_\Phi(x_h)\boldsymbol{e}_3, \\
&-\Dl_h g^i_\Phi(x_h)=C_i,\q\q\q\text{in }\q \Sigma_i,\\
&\f{\p g^i_\Phi}{\p\bl{\bar{n}}}=-\al g^i_\Phi,\q\q\q\hskip .7cm\text{on }\q \p\Sigma_i,\\
&\int_{\Sigma_i}g^i_\Phi(x_h)dx_h=\Phi, \\
\end{aligned}
\right.
\]
where the constant $C_i$ is uniquely related to $\Phi$, while $\bl{\bar{n}}$ is the unit outer normal vector on $\p\Sigma_i$. We can see that $\bl{g}^i_{\Phi}$ is a solution of \eqref{NS} with the Navier-slip boundary \eqref{NBC} in $\Sigma_i\times\bR$.

The main objective of this paper is to study the solvability of the following generalized Leray's problem: For a given flux $\Phi$, to find a pair $(\bl{u},p)$ such that
\be\label{GL1}
\lt\{
\bali
&\bl{u}\cdot\na \bl{u}+\na p-\Dl \bl{u}=0,\q \na\cdot \bl{u}=0,\q \hskip .1cm\text{in}\q \mathcal{D}, \\
&2(\mathbb{S}\bl{u}\cdot\boldsymbol{n})_{\mathrm{tan}}+\alpha \bl{u}_{\mathrm{tan}}=0,\q \bl{u}\cdot\boldsymbol{n}=0, \q \text{on}\q \p\mathcal{D},\\
\eali
\rt.
\ee
with
\be\label{GL2}
\int_{\Sigma_i}u_3(x_h, x_3) d x_h=\Phi,\q\text{for}\q |x_3|>Z/2,
\ee
and
\be\label{GL3}
\bl{u}\rightarrow \boldsymbol{g}^i_{\Phi},\q \text{as}\q |x|\rightarrow \i \ \text{ in }\ \mD_i.
\ee
To prove the existence of the above generalized Leray's problem, we first introduce a weak formulation. Multiplying \eqref{GL1}$_1$ with  $\bl{\varphi}\in\bl{X}$ and integration by parts, by using the boundary condition \eqref{GL1}$_2$, we can obtain
\be\label{weaksd}
\begin{split}
&2\int_{\mD}\mathbb{S}\bl{u}:\mathbb{S}\bl{\varphi} dx+\al\int_{\p\mD}\bl{u}_{\mathrm{tan}}\cdot\bl{\varphi}_{\mathrm{tan}}dS+\int_{\mD}\bl{u}\cdot\nabla \bl{\varphi} \cdot \bl{u} dx=0,\q \text{for all }\ \bl{\varphi}\in \bl{X}.
\end{split}
\ee
Now we define the weak solution of the generalized Leary's problem:

\begin{definition}\label{weaksd1}
A vector $\bl{u}:\mD\to \bR^3$ is a called a {\it weak} solution of the generalized Leray problem \eqref{GL1} to \eqref{GL3} if and only if
\begin{itemize}
\setlength{\itemsep}{-2 pt}
\item[(i).] $\bl{u}\in \mathcal{H}_{\sigma,{\mathrm{loc}}}(\overline{\mD})$;

\item[(ii).] $\bl{u}$ satisfies \eqref{weaksd};

\item[(iii).] $\bl{u}$ satisfies \eqref{GL2} in the trace sense;

\item[(iv).] $\bl{u}-\bl{g}^i_{\Phi}\in {H}^1(\mD_i)$, for $i=L,R$.
\end{itemize}
\end{definition}

\qed

\begin{remark}
The weak solution also satisfies a generalized version of \eqref{GL3}. Actually, it follows from the trace inequality (\cite[Theorem II.4.1]{Galdi2011}) that for any $x_3>Z$:
\[
\bali
&\int_{\Sigma_R}|\bl{u}(x_h,x_3)-\bl{g}^R_{\Phi}(x_h)|^2dx_h\leq C\|\bl{u}-\bl{g}^R_{\Phi}\|^2_{H^1(\Sigma_R\times(x_3,+\i))},
\eali
\]
where the constant $C$ is independent of $x_3$. This implies that
\[
\int_{\Sigma_R}|\bl{u}(x_h,x_3)-\bl{g}^R_{\Phi}(x)|^2dx_h \to0,\quad\text{as $x_3\to\infty$}.
\]
The case of $x_3<-Z$ is similar.
\end{remark}

\qed

The following result shows that for each weak solution we can associate a corresponding pressure field. See the proof in Section \ref{SEC332} below.

\begin{lemma}\label{pressure}
Let $\bl{u}$ be a weak solution to the generalized Leray's problem defined above. Then there exists a scalar function $p\in L^2_{{\mathrm{loc}}}(\overline{\mD})$ such that
\[
\int_{\mD}\na\bl{u}:\na\bl{\psi} dx+\int_{\mD}\bl{u}\cdot\nabla \bl{u}\cdot \bl{\psi}dx=\int_{\mD} p\na\cdot\bl{\psi} dx
\]
holds for any $\bl{\psi}\in C^\i_c(\mD; \bR^3)$.
\end{lemma}

\qed

\subsection{Main results}
\q\ Now we are ready to state the main theorems of this paper. The first one is the existence of weak solutions, the second one addresses uniqueness of the weak solution and the third one concerns regularity and decay estimates of the weak solution:
\begin{theorem}\label{Ext}
Assume that $\mD$ is the aforementioned smoothness domain in \eqref{cylinder}. Then there exists a positive constant $\Phi_0$, depending only on $\al$ and $\mD$, such that for any $\Phi\leq \Phi_0$, the generalized Leray's problem \eqref{GL1}-\eqref{GL3} has a weak solution $(\bl{u},p)\in \mathcal{H}^1_{\sigma,{\mathrm{loc}}}(\overline{\mD})\times L^2_{{\mathrm{loc}}}(\overline{\mD})$ satisfying
\be\label{sesti}
\sum_{i=L,R}\|\bl{u}-\boldsymbol{g}^i_\Phi\|_{H^1(\mD_i)}\leq C_{\al,\mD}\Phi,
\ee
where $C_{\al,\mD}$ depends only on $\al$ and $\mD$.
\end{theorem}

\qed

Theorem \ref{Ext} only states the existence of a weak solution, which has finite energy property after subtracting the background Poiseuille flows. However, uniqueness of the solution is not shown. Actually, estimate \eqref{sesti} is far enough for us to deduce the following theorem of uniqueness, which allows the energy of the weak solution in $\mD_\zeta$ satisfies a $3/2-order$ growth with respect to $\zeta$. The achievement of this relaxation is due to the application of a key Lemma in \cite{Lady-Sol1980}, which will be presented in Section \ref{SEC2}.

\begin{theorem}\label{thqunique}
Let $(\bl{u},p)$ be a weak solution to \eqref{GL1}--\eqref{GL3}. Suppose for any $\zeta>Z$,
\be\label{GROWC}
\|\na \bl{u}\|_{L^2(\mD_\zeta)}=o\lt(\zeta^{3/2}\rt),
\ee
Then the weak solution is unique \footnote{ The pressure $p$ is unique up to subtracting an arbitrary constant.} provided the flux $\Phi$ is sufficiently small.
\end{theorem}

\qed

The following Theorem gives the smoothness and the asymptotic behavior of $(\bl{u},p)$, which decays exponentially to the Poiseuille flow $\boldsymbol{g}^i_{\Phi}$ at each outlet $\mD_i$ as $x_3$ tends  to infinity.

%

\begin{theorem}\label{THM16}
Let $\bl{u}$ be the weak solution in Theorem \ref{Ext} and $p$ is the corresponding pressure. Then
\[
(\bl{u},p)\in C^\i(\overline{\mD})
\]
such that: For any  $m=|\beta|\geq 0$,
\be\label{HODEST}
\sum_{i=L,R}\|\na^\beta(\bl{u}-\boldsymbol{g}^i_\Phi)\|_{L^2(\mD_i)}+\|\na^\beta \bl{u}\|_{L^2(\mD_M)}\leq C_{m,\al,\mD}\Phi.
\ee
Meanwhile, the following pointwise decay estimates  hold:
\be\label{pointdecay}
\begin{split}
|\na^\beta(\bl{u}-\boldsymbol{g}^L_{\Phi})(x)|&\leq C_{m,\al,\mD}\Phi \exp\left\{-\sigma_{m,\al,\mD}|x_3|\right\},\q\text{for all}\q x_3<-Z-1;\\
|\na^\beta(\bl{u}-\boldsymbol{g}^R_{\Phi})(x)|&\leq C_{m,\al,\mD}\Phi \exp\left\{-\sigma_{m,\al,\mD}|x_3|\right\},\q\text{for all}\q x_3>Z+1.
\end{split}
\ee
Here $C_{m,\al,\mD}$ and $\sigma_{m,\al,\mD}$ are positive constants depending on $m$, $\al$ and $\mD$.
\end{theorem}

\qed

\begin{remark}
For the pressure $p$ which is generated in Lemma \ref{pressure}, there exists two constants $C_{P,L},\,C_{P,R}>0$ (See \eqref{POS}) , and a smooth cut-off function $\eta$ with
\[
\eta(x_3)=\left\{
\begin{aligned}
&1,\q\text{for}\q x_3>Z;\\
&0,\q\text{for}\q x_3<Z/2,\\
\end{aligned}
\right.
\]
(which is given in \eqref{ETA}) such that: For any  $m=|\beta|\geq 0$,
\[
\left\|\na^\beta\na\left(p+\f{\Phi\int_{-\i}^{x_3}\eta(s)ds}{C_{P,R}}-\f{\Phi\int_{-\i}^{-x_3}\eta(s)ds}{C_{P,L}}\right) \right\|_{L^2(\mD)}\leq C_{m,\al,\mD}\Phi.
\]
Meanwhile, the following pointwise decay estimate holds: for all $|x_3|>Z+1$,
\[
\left|\na^{\beta}\na \left(p+\f{\Phi\int_{-\i}^{x_3}\eta(s)ds}{C_{P,R}}-\f{\Phi\int_{-\i}^{-x_3}\eta(s)ds}{C_{P,L}}\right)(x)\right|\leq C_{m,\al,\mD}\Phi \exp\left\{-\sigma_{m,\al,\mD}|x_3|\right\},
\]
where $C_{m,\al,\mD}$ and $\sigma_{m,\al,\mD}$ are positive constants depending on $m$, $\al$ and $\mD$. The subtracted term
\[
p_{\bl{g}}:=-\f{\Phi\int_{-\i}^{x_3}\eta(s)ds}{C_{P,R}}+\f{\Phi\int_{-\i}^{-x_3}\eta(s)ds}{C_{P,L}}
\]
is set to balance the pressure of the Poiseuille flows.
\end{remark}
\qed

\subsection{Main difficulties, strategies and outline of the proof}
\subsubsection*{Difficulties}
\q\ Compared with the no-slip boundary condition, the main difficulties of the problem with Navier-slip boundary condition lie in the following:
\begin{itemize}
\setlength{\itemsep}{-3 pt}
\item[(i).] The absence of the Korn-type inequality ($L^2$ norm equivalence between $\na\bl{v}$ and $\mathbb{S}\bl{v}$) on $\mD$ with noncompact boundary;

\item[(ii).] For a given flux, construction of a smooth solenoidal vector field $\bl{a}$, satisfying the Navier-slip boundary condition, and equalling to the Poiseuille flow at a large distance;

\item[(iii).] Achieving Poincar\'e-type inequalities under the Navier-slip boundary condition;

\item[(iv).] Derivation of the global $H^2$ estimate of the $H^1$-weak solution.
\end{itemize}

\subsubsection*{Strategies}
\q\ In order to overcome difficulties listed above, our main strategies are as follows:

\begin{itemize}
\setlength{\itemsep}{-3 pt}
\item[(i).] During the proof of both existence and uniqueness, an $H^1$ estimate of the solution is required. Owing to the boundary condition, we only have the $L^2$ estimate of stress tensor $\mathbb{S}v$. However, the Korn-type inequality is not applicable for our domain considered here with a noncompact boundary. Fortunately, the energy estimate of stress tensor $\mathbb{S}v$ will produce a boundary integration with a good sign, which can be used to control the bad terms coming from the energy estimate of the gradient of the velocity. At last, by combining the uniform energy estimates of the stress tensor with the $L^2$-estimate of the gradient of the velocity, we can achieve the $H^1$ estimate of $\bl{v}$.

\item[(ii).] Our main idea in constructing $\bl{a}$ is to smoothly connect Poiseuille flows in $\bl{g}_\Phi^L$ and $\bl{g}_\Phi^R$ with a compact supported divergence-free vector $(0,0,h(x_h))$ in $\mD_M$. In the intermediate parts, we glue them by solving a 2D divergence equation in the cross section with the 2D Navier-slip boundary condition.

\item[(iii).] For the no-slip boundary condition, the Poincar\'e inequality can be applied directly. However, in the case of the Navier-slip boundary condition, the Poincar\'e inequality is not obvious in both straight pipes $\mD_L$, $\mD_R$ and the truncated pipe $\mD_\zeta$. To handle the case in $\mD_L$ or $\mD_R$, we divide a vector-valued function into $x_h$-direction part and $x_3$-direction part. The first part follows from a 2D Payne's identity \eqref{Payne} and impermeable boundary condition, while the second part is achieved by subtracting the constant flux so that $v_3$ has zero mean value in any cross section of $\mD_L$ or $\mD_R$. See Lemma \ref{P2}. Based on the result of the straight pipe, we derive the Poincar\'e inequality in $\mD_\zeta$ by the trace theorem and a 3D Payne's identity \eqref{Payne3D}. See Lemma \ref{TORPOIN}.

\item[(iv).] Our idea of obtaining the global $H^2$ estimate is to decompose $\mD$ into a series of bounded smooth domains $\tilde{\mfD}_k$ which only have three shapes, so that the related estimate constant in $\tilde{\mfD}_k$ could be uniform with $k$. In each $\tilde{\mfD}_k$, we establish the $H^2$ estimate of the solution via the known conclusions for the linear Stokes system with the Navier-slip boundary condition in \cite{AACG:2021JDE}. Then we achieve the global $H^2$ estimate by summarizing those estimates in $\tilde{\mfD}_k$.
\end{itemize}

\subsubsection*{Outline of the proof}
\q\ The existence of the solution will be given in Section \ref{SECE}. First, the Poiseuille flows in $\mD_i$ ($i=L,R$), with their fluxes being $\Phi$ and satisfying the Navier-slip boundary condition will be constructed. Then a smooth divergence-free vector field in $\mD$, subject to the Navier-slip boundary condition and equaling to the Poiseuille flows at the far left and far right, will be introduced. In this way, we can reduce the existence problem to a related one that the solution approaches zero at spacial infinity. Then this problem can be handled by the standard Galerkin method.

The proof of the uniqueness is derived in Section \ref{sec2}. The main idea is applying Lemma \ref{LEM2.3}, which was originally announced in reference \cite{Lady-Sol1980} as far as the authors know. If $(\bl{u},p)$ and $(\tilde{\bl{u}}, \tilde{p})$ are two distinct solutions, we denote the energy integral in terms of $\bl{w}:=\tilde{\bl{u}}-\bl{u}$ as follows:
\bes
Y(K):=\int_{K-1}^K\int_{\mathcal{D}_\zeta}|\nabla \bl{w}|^2dxd\zeta.
\ees
An ordinary differential inequality of $Y(K)$, which satisfies the assumption in Lemma \ref{LEM2.3}, will be derived. The derivation of this inequality involves in a series of estimates, two terms of which are especially different from the previous literature, $i.e.$,  $\int_{\O_\zeta} p v_3dx$ and $\int_{\mD_\zeta}\bl{v}\cdot\Dl \bl{v}dx$. The estimate of the first term involves in an application of partial Poinc\'{a}re inequality in Lemma \ref{P1} and a divergence-gradient operator estimate in Lemma \ref{LEM2.1}. The estimate of the second term is derived by combining the $L^2$ norm of both the stress tensor and the gradient of the velocity. At last, the vanishing of $Y(K)$ will be proved, which indicates the uniqueness of the solution.

Proofs of the smoothness and exponential decay of the solution to the Poiseuille flows at large spacial distance are given in Section \ref{SEC5}. By using the ``decomposing-summarizing" technique in the Strategies (iv), a $W^{1,3}$ and then $H^2$ global estimates will be derived. Using the bootstrapping argument, higher-order global estimates will follow. For the exponential decay estimates, here goes the idea: By defining the following energy as, for $\zeta>Z+1$,
\bes
\mathcal{G}(\zeta):=\int_{\Sigma_L\times(-\i,-\zeta)}|\na(\bl{u}-\bl{g}^{L}_{\Phi})|^2dx+\int_{\Sigma_R\times(\zeta,+\i)}|\na(\bl{u}-\bl{g}^{R}_{\Phi})|^2dx,
\ees
then we can derive a first order ordinary differential inequality, which will result in the exponential decay of $\mathcal{G}(\zeta)$. Finally, higher-order estimates of the solution in $\mD-\mD_\zeta$, Sobolev embedding, and the exponential decay of $\mathcal{G}(\zeta)$ will validate the pointwise decay of the solution in \eqref{pointdecay}.

\subsection{Related works}
\q\ Before ending the introduction, we review some works related to the solvability of the Leray's problem in \eqref{GL1}, \eqref{GL2} and \eqref{GL3}. The original Leray's problem concerns the existence, uniqueness, regularity and asymptotic behavior of System \eqref{GL1}--\eqref{GL3} with no-slip boundary condition (corresponding to $\al=+\i$ in \eqref{GL1}$_2$). See the description in Ladyzhenskaya \cite[p. 77]{Ladyzhenskaya:1959UMN} and \cite[p. 551]{Ladyzhenskaya:1959SPD}. Amick \cite{Amick:1977ASN, Amick:1978NATMA} contributed the first remarkable work on the solvability of the Leray's problem with small flux and the no-slip boundary condition, which reduced the solvability problem to the resolution of a variational problem. However, uniqueness is left open. A detailed analysis of the existence, uniqueness and asymptotic behavior of small-flux solutions is given by Ladyzhenskaya and Solonnikov \cite{Lady-Sol1980}. More details on well-posedness, decay and far-field asymptotic analysis of solutions for the Leray's problem with the no-slip boundary and related topics can be found in \cite{AP:1989SIAM, HW:1978SIAM, Pileckas:2002MB} and references therein. Readers can trace to \cite[Chapter XIII]{Galdi2011} for a systematic review and study of the Leray's problem with the no-slip boundary. Recently the third author of the present paper and Yin \cite{YY:2018SIAM} studied the well-posedness of weak solutions to the steady non-Newtonian fluids in pipe-like domains. Wang-Xie in \cite{WX:2019ARXIV,Wang-Xie2022ARMA} studied the existence, uniqueness and uniform structural stability of Poiseuille flows for the 3D axially symmetric inhomogeneous (a force term appearing on the right hand of the equations) Navier-Stokes equations in the 3D pipe.

Compared to the no-slip boundary condition, the Leray's problem with the Navier-slip boundary condition, which also has different physical interpretations and mathematical properties, seems to be much more complicated. Literature \cite{Mucha2003, Mucha:2003STUDMATH, Konie:2006COLLMATH} studied the solvability of the steady Navier-Stokes equations with the perfect Navier-slip condition ($\al=0$), where they employed a constant vector field as its asymptotic profile at the spatial infinity. Only the existence, regularity, and asymptotic behavior of weak solutions were addressed there. The uniqueness was left open, and the asymptotic behavior at far fields was not given there. The problem raised there could be recognized as the Leray's problem with the complete Navier-slip boundary condition in a two- dimensional strip, where the asymptotic profile is a constant vector. Our problem raised in \eqref{GL1}--\eqref{GL3} is a perfect extension of the original Leray's problem with the no-slip boundary. The background Poiseuille flows considered here tend to the Leray's Poiseuille flows with the no-slip boundary condition as $\al\to +\i$. As far as the authors know, there is little literature concerning the solvability of the generalized Leray's problem \eqref{GL1}--\eqref{GL3}, which settles the well-posedness issue on the steady Navier-Stokes equations subject to the Navier-slip boundary in an unbounded domain with an unbounded boundary. While for the well-posedness of solutions to the steady Navier-Stokes equations with the Navier-slip boundary in bounded domains, there already have many works that we can refer to, see \cite{Amrouche2014,Ghosh:2018PHD,AACG:2021JDE} and references therein. Recently, Wang and Xie \cite{WX:2021ARXIV2} gave the uniqueness and uniform structural stability of Poiseuille flows in an infinitely long pipe with the Navier-slip boundary condition for the inhomogeneous axially symmetric Navier-Stokes equations. The first and second authors of the present paper gave the characterization of bounded smooth solutions for the axially symmetric Navier-Stokes equations with the perfect Navier-slip boundary condition (corresponding to $\al=0$ in \eqref{NBC}$_1$) in the infinitely long cylinder \cite{LP:2021ARXIV}.

This paper is arranged as follows: Section \ref{SEC2} contains preliminary work of the proof, in which the Navier-slip boundary condition will be written under the ``natural" moving frame of $\p\mD$, and some useful lemmas will be presented. Section \ref{SECE} is devoted to the proof of existence results. In Section \ref{sec2}, we will finish the proof of the uniqueness of the solution. Finally, we focus on the higher-order regularity and exponential decay properties of the solution in Section \ref{SEC5}.

At last, we emphasize that the domain considered in this paper is a 3-D  distorted pipe, while, if the domain is a two dimensional distorted strip, similar results as stated in Theorem \ref{Ext}, \ref{thqunique} and \ref{THM16} will also be obtained. More precisely, if we consider the Navier-Stokes equation \eqref{NS} with the Navier-slip boundary condition \eqref{NBC} in the strip $[0,1]\times\bR$, the following two dimensional Poiseullie flow will be obtained by a direct calculation
\be\label{2dpoi}
\bl{g}_{\Phi}=\lt(0, \f{6\al\Phi}{6+\al}\lt(-x^2_1+x_1\rt)+\f{6\Phi}{6+\al}\rt).
\ee
After a compact perturbation of the domain $\bR\times[0,1]$, the existence, uniqueness, regularity of the solutions, which approaches to $\bl{g}_{\Phi}$ in \eqref{2dpoi} at spacial infinity will be presented in our forthcoming paper, where the flux at the cross section $\Phi$ can be relatively large.

\section{Preliminary}\label{SEC2}
\subsection{Reformulation of the boundary condition in the local orthogonal curvilinear coordinates}
\q\ First, we rewrite the boundary condition \eqref{NBC} in the locally moving coordinate framework.

Regarding the smoothness of the pipe $\mD$, for any given point on $\p\mD$, we define $(\gamma_1,\gamma_2,\gamma_3)$ a system of orthogonal curvilinear coordinates in $U\subset\mathbb{R}^3$, where $U$ is a neighborhood of the aforementioned point. The surface $\gamma_3=0$ represents a portion of surface $\mD$, and surfaces $\gamma_3=constant$ are parallel to this portion with $\gamma_3$ increases towards the outside of $\mD$. On each surface $\gamma_3={constant}$, two family of curves, the  $\gamma_1-{curve}$ and $\gamma_2-{curve}$, are lines of curvature of the surface. Their unit tangent vectors $\bl{\tau_1}$, $\bl{\tau_2}$ and the normal vector $\boldsymbol{n}$ form an orthogonal basis at each point of the neighborhood $U$, with the Lam\'e coefficients $H_1,H_2,H_3>0$ such that
\[
\left\{
\begin{aligned}
&\p_{\gamma_i}x=H_i\bl{\tau_i},\q\text{for }i=1,2;\\
&\p_{\gamma_3}x=H_3\bl{n}.\\
\end{aligned}
\right.
\]
\begin{figure}[H]
\centering
\includegraphics[scale=0.5]{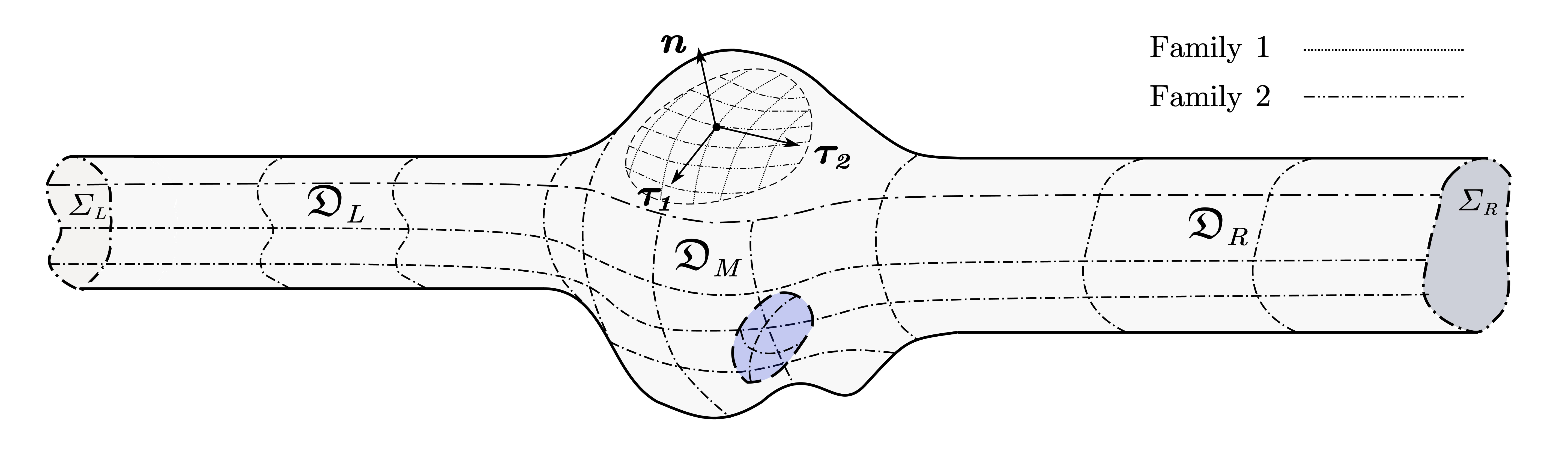}
\caption{The local orthogonal curvilinear coordinates on $\p\mD$}
\end{figure}

Under this (local) curvilinear coordinates, one can write
\[
\bl{u}=u_{{\tau_1}}\bl{\tau_1}+u_{\tau_2}\bl{\tau_2}+u_n\bl{n}.
\]
Then \eqref{NBC} enjoys the following simplified expression:
\be\label{NBCM}
\left\{
\begin{aligned}
&\p_{\boldsymbol{n}}u_{\tau_1}=\left(\kappa_1(x)-\al\right)u_{\tau_1},\\
&\p_{\boldsymbol{n}}u_{\tau_2}=\left(\kappa_2(x)-\al\right)u_{\tau_2},\\
&u_n=0,\\
\end{aligned}
\right.\q\text{on}\ U\cap\p\mD.
\ee
See \cite[Proposition 2.1 and Corollary 2.2]{Watanabe2003}. Here for $i=1,2$,
\be\label{curvature}
\kappa_i(x)=-\f{\bl{n}}{H_i}\cdot \f{\p\bl{\tau}_i}{\p\g_i}
\ee
are the principal curvatures of $\p\mD$ corresponding to $\gamma_i-{curve}$, respectively.

Noting that $\p\mD$ can be divided into two parts. The part $\p\mD\cap\p{\mD_M}$ is a compact manifold, thus one can deduce the uniform boundedness of $\kappa_i$  ($i=1,2$) there. The other part $\p\mD\cap\p\left(\mD_L\cup\mD_R\right)$ is a combination of two semi-infinite smooth straight pipes, whose curvature depends only on the scalar curvature of smooth Jordan curves $\p\Sigma_L$ and $\p\Sigma_R$, which is clearly uniformly bounded. More details on the locally natural moving frame on $\p\mD\cap(\p\mD_R\cup\p\mD_L)$ is shown in the following remark:

\begin{remark}\label{RMK2.2}
In the ``straight" part of the pipe (i.e. $\mD_L\cup\mD_R$), one can always choose $\bl{\tau_2}=\bl{e_3}$, which is constant vector. Meanwhile, $\bl{\tau_1}$ and $\bl{n}$ are unit tangent vector and unit outer normal vector of the cross section $\Sigma$, respectively. See below
\begin{figure}[H]
\centering
\includegraphics[scale=0.6]{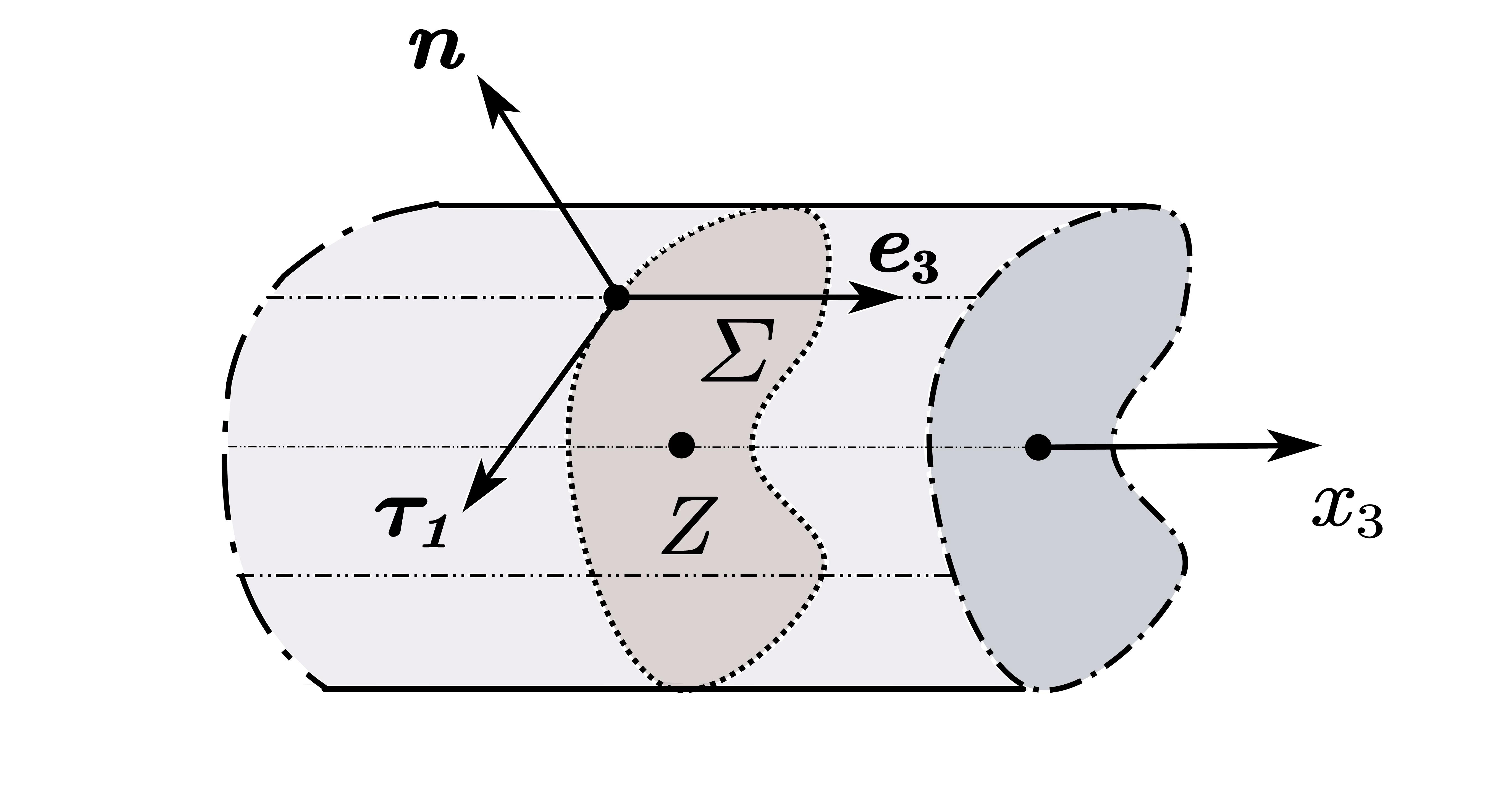}
\caption{The orthogonal curvilinear coordinates on $\p\mD\cap(\p\mD_R\cup\p\mD_L)$}
\end{figure}
In this case, we find $\boldsymbol{\tau_1}$ and $\boldsymbol{n}$, which are perpendicular to the $x_3$-axis, depend only on $(x_1,x_2)$. And from \eqref{curvature}, the principal curvature $\kappa_2$ is identical zero. By writing
\[
\bl{u}=u_{{\tau_1}}\bl{\tau_1}+u_3\bl{e_3}+u_n\bl{n},
\]
\eqref{NBCM} is simplified to
\be\label{NBCM1}
\left\{
\begin{aligned}
&\p_{\boldsymbol{n}}u_{\tau_1}=\left(\kappa_1(x)-\al\right)u_{\tau_1},\\
&\p_{\boldsymbol{n}}u_3=-\al u_3,\\
&u_n=0,\\
\end{aligned}
\right.\q\q\q\text{on}\q\p\mD\cap(\p\mD_R\cup\p\mD_L).
\ee
Since $\p\Sigma_i$ ($i=L$ or $R$) is smooth and compact, the principal curvature $\kappa_1$ must be uniformly bounded on $\p\mD\cap(\p\mD_R\cup\p\mD_L)$.
\end{remark}

\qed

Thus one concludes the following result at the end of this subsection:
\begin{proposition}\label{Prop22}
The principal curvature $\kappa_i$ ($i=1$ or $2$) is uniformly bounded on $\p\mD$.
\end{proposition}

\qed

\subsection{Useful lemmas}
\q\ In this part, we give some useful lemmas which will be frequently used throughout the rest of this paper. Lemma \ref{P0}--Lemma \ref{TORPOIN} concern on Poincar\'e inequalities of the solution $\bl{u}$ in a part of the straight pipe $\Sigma\times\mathbb{R}$ and truncated finite pipe $\mD_\zeta$ with only the impermeable condition \eqref{NBC}$_2$. Lemma \ref{LEM2.3} is introduced to show the uniqueness of the solution. Lemma \ref{LEM52} and Lemma \ref{LEM59} are regularity results of linear Stokes systems on a bounded domain, which will be applied in the bootstrapping argument in Section \ref{SEC5}.

In the standard Euclidean coordinates framework $\boldsymbol{e_1},\ \boldsymbol{e_2}$ and $\boldsymbol{e_3}$, Let $\bl{u}=u_1\boldsymbol{e_1}+u_2\boldsymbol{e_2}+u_3\boldsymbol{e_3}$. $\Sigma\subset\bR^2$ is a smooth bounded domain in $x_h$ directions and $I\subset\bR$ is an (infinite or finite) interval in $x_3$ direction. We set $\boldsymbol{n}=(n_1,n_2,0)$ is the unit outer normal vector on $\p\Sigma\times I$, where $\bar{\boldsymbol{n}}=(n_1,n_2)$ is the unit outer normal vector on $\p\Sigma$.

\begin{lemma}\label{P0}
Let $\Sigma\subset\mathbb{R}^2$ be a compact domain with $C^1$ boundary, and $\bl{f}=f_1\boldsymbol{e_1}+f_2\boldsymbol{e_2}$ be a two dimensional vector function with components in $H^1(\Sigma)$, and $\bl{f}\cdot\bar{\boldsymbol{n}}=0$ on $\partial \Sigma$, where $\bar{\boldsymbol{n}}$ is the unit outer normal of $\p \Sigma$. Then the following Poincar\'e inequality holds
\be\label{parpoin}
\|\bl{f}\|_{L^2(\Sigma)}\leq C_\Sigma\|\nabla_h \bl{f}\|_{L^2(\Sigma)},
\ee
where $\na_h=(\p_{x_1},\p_{x_2})$ is the gradient operator on $x_1$ and $x_2$ direction and $C_\Sigma=Cd_{\Sigma}$, where $C$ is an absolute constant and $d_{\Sigma}:=\max_{x_h,y_h\in\Sigma}\{|x_h-y_h|\}$ is the diameter of $\Sigma$.
\end{lemma}
\pf  See the hint in Galdi \cite[Page 71, Exercise II.5.6]{Galdi2011}. Here, we give the proof for completeness. First we choose a fixed point $x_0=(x_{0,1}, x_{0,2})\in \Sigma$, then it is easy to check the following equality:
\be\label{Payne}
\sum^2_{i,j=1}\lt[\p_{x_i}(f_i (x_j-x_{0,j}) f_j)-\p_{x_i} f_i(x_j-x_{0,j}) f_j-|\bl{f}|^2-f_i(x_j-x_{0,j})\p_{x_i} f_j\rt]=0.
\ee
Integrating the above equality on $\Sigma$, one deduces
\be\label{ZJ1}
\bali
\int_\Sigma |\bl{f}|^2dx_h=&\underbrace{\sum^2_{i,j=1}\int_\Sigma\p_i(f_i (x_j-x_{0,j}) f_j)dx_h}_{I_1}-\sum^2_{i,j=1}\int_\Sigma\p_i f_i(x_j-x_{0,j}) f_jdx_h\\
                      &-\sum^2_{i,j=1}\int_\Sigma f_i(x_j-x_{0,j})\p_i f_jdx_h.\\
\eali
\ee
Using the divergence theorem and the boundary condition $\bl{f}\cdot\bar{\boldsymbol{n}}=0$, we can obtain
\[
I_1=\sum^2_{j=1} \int_{\p\Sigma} \bar{\boldsymbol{n}}\cdot \bl{f}(x_j-x_{0,j}) f_jdS=0.
\]
Thus by \eqref{ZJ1} and the Cauchy-Schwarz inequality, we arrive at
\[
\bali
\int_\Sigma |\bl{f}|^2dx_h\leq \f{1}{2}\int_\Sigma |\bl{f}|^2dx_h+Cd^2_\Sigma\int_\Sigma |\na_h \bl{f}|^2dx_h,
\eali
\]
which indicates \eqref{parpoin}.

\qed

If we choose  $\bl{f}=u_1\boldsymbol{e_1}+u_2\boldsymbol{e_2}$ as in Lemma \ref{P0}, we deduce the following Lemma:

\begin{lemma}[Partial Poincar\'e inequality in a straight pipe]\label{P1}
Let $\bl{u}=u_1\boldsymbol{e_1}+u_2\boldsymbol{e_2}+u_3\boldsymbol{e_3}$ be a $H^1$ vector field in $\Sigma\times I$. If $\bl{u}$ satisfies the boundary condition $\bl{u}\cdot\boldsymbol{n}=0$, then the following Poincar\'e inequality holds
\be\label{POV1}
\|u_1\boldsymbol{e_1}+u_2\boldsymbol{e_2}\|_{L^2(\Sigma\times I)}\leq C\|\nabla_h\left(u_1\boldsymbol{e_1}+u_2\boldsymbol{e_2}\right)\|_{L^2(\Sigma\times I)},
\ee
where $C$ is an absolute positive constant.
\end{lemma}

\pf For any $x=(x_h,x_3)\in\Sigma\times I$, by the impermeable condition $\bl{u}\cdot\boldsymbol{n}=0$, one sees
\[
\left(u_1\boldsymbol{e_1}+u_2\boldsymbol{e_2}\right)(x_h,x_3)\cdot\bar{\boldsymbol{n}}=\bl{u}\cdot\boldsymbol{n}=0,\q\text{for any } x_h\in\p\Sigma.
\]
Then Lemma \ref{P0} indicates
\[
\|(u_1\boldsymbol{e_1}+u_2\boldsymbol{e_2})(\cdot,x_3)\|^2_{L^2(\Sigma)}\leq C^2\|\nabla_h(u_1\boldsymbol{e_1}+u_2\boldsymbol{e_2})(\cdot,x_3)\|^2_{L^2(\Sigma)}.
\]
Integrating with $x_3$ on $I$ on both sides, respectively, one concludes \eqref{POV1}.

\qed

We notice that the Poincar\'e inequality could not hold for $u_3$, due to the existence of parallel flow (See e.g. \eqref{HP1} below). Nevertheless, the mean value of $u_3$ through the cross section $\Sigma$ is conserved for $x_3\in I$ if $\bl{u}$ is a divergence-free vector field. The reason is: Denoting the flux flowing across $\Sigma$ by
\[
\Phi(x_3)=\int_{\Sigma}\bl{u}(x_h,x_3)\cdot\boldsymbol{e_3}dx_h,
\]
then by applying the divergence-free property of $\bl{u}$ and the impermeable condition $\bl{u}\cdot\boldsymbol{n}=0$, we have
\[
\begin{split}
\f{d}{dx_3}\Phi(x_3)&=\int_{\Sigma}\p_{x_3}u_3(x_h,x_3)dx_h\\
&=-\int_{\Sigma}(\p_{x_1}u_1+\p_{x_2}u_2)(x_h,x_3)dx_h\\
&=-\int_{\p\Sigma}(\bl{u}\cdot\boldsymbol{n})(x_h,x_3)dS(x_h)=0.
\end{split}
\]
This implies the constancy of the flux $\Phi$. Then we have the following Lemma:
\begin{lemma}\label{P2}
Let $\bl{u}=u_1\boldsymbol{e_1}+u_2\boldsymbol{e_2}+u_3\boldsymbol{e_3}$ be a $H^1$ vector field in $\Sigma\times I$, which is divergence-free and satisfies the boundary condition $\bl{u}\cdot\boldsymbol{n}=0$. Set $g(x_h)\in H^1(\Sigma)$ satisfying
\[
\int_{\Sigma}g(x_h)dx_h=\Phi.
\]
and define $\bl{v}:=\bl{u}-g(x_h)\bl{e_3}$, then we have the following

\be\label{poinfull}
\|\bl{v}\|_{L^2(\Sigma\times I)}\leq C\|\nabla_h \bl{v}\|_{L^2(\Sigma\times I)}.
\ee
where $C$ is an absolute positive constant.
\end{lemma}

\pf
Notice that $v_3$ has a vanishing mean value on the cross section $\Sigma$, and therefore it enjoys the 2D Poincar\'e inequality:
\be\label{Ex}
\bali
\int_{\Sigma}|v_3(x_h,x_3)|^2dx_h=&\int_{\Sigma}\left|v_3(x_h,x_3)-\f{1}{|\Sigma|}\int_{\Sigma} v_3(x_h,x_3)dx_h\right|^2dx_h\leq C^2\int_{\Sigma}|\na_h v_3(x_h,x_3)|^2dx_h.
\eali
\ee
This indicates the 3D Poincar\'e inequality
\[
\|v_3\|_{L^2(\Sigma\times I)}\leq C\|\nabla_hv_3\|_{L^2(\Sigma\times I)}
\]
by integrating \eqref{Ex} with $x_3$ on the interval $I$. Combining the results in Lemma \ref{P1}, we can obtain \eqref{poinfull}.

\qed


Based on Lemma \ref{P2}, one has the following Poincar\'e-type inequality in the truncated distorted pipe $\mD_\zeta=\{x\in\mD:\,-\zeta\leq x_3\leq \zeta\}$.
\begin{lemma}\label{TORPOIN}
Given $\zeta\geq Z$, let $\bl{w}=(w_1,w_2,w_3)\in H^1(\mD_\zeta)$ with zero flux in $\mD_\zeta$, that is
\[
\int_{\Sigma_R}\bl{w}(x_h,Z/2)\cdot\bl{e_3}dx_h=0.
\]
If we suppose $\bl{w}\cdot\bl{n}\equiv 0$ on $\p\mD\cap\p\mD_\zeta$, where $\bl{n}$ is the unit outer normal vector on $\p\mD$, then the following Poincar\'e inequality holds:
\be\label{TORPIPEPOIN}
\|\bl{w}\|_{L^2(\mD_\zeta)}\leq C_\mD\|\na\bl{w}\|_{L^2(\mD_\zeta)}.
\ee
Here $C_\mD>0$ is a constant which is uniform with $\zeta$.
\end{lemma}
\pf Integrating the following identity on $\mD_M=\{x\in\mD :\,-Z/2\leq x_3\leq Z/2\}$
\be\label{Payne3D}
\sum^3_{i,j=1}\lt[\p_{x_i}(w_i x_j w_j)-\p_{x_i} w_ix_jw_j-|\bl{w}|^2-w_ix_j\p_{x_i} w_j\rt]=0,
\ee
and using a similar approach as we estimate terms on the right hand side of \eqref{ZJ1}, one derives
\[
\begin{split}
\int_{\mD_M}|\bl{w}|^2dx\leq& \f{1}{2}\int_{\mD_M}|\bl{w}|^2dx+C_\mD\int_{\mD_M}|\na\bl{w}|^2dx\\
&+\left|\int_{\Sigma_L}\left(w_3(x\cdot\bl{w})\right)(x_h,-Z/2)dS\right|+\left|\int_{\Sigma_R}\left(w_3(x\cdot\bl{w})\right)(x_h,Z/2)dS\right|,
\end{split}
\]
which indicates
\be\label{PCC1}
\int_{\mD_M}|\bl{w}|^2dx\leq C_\mD\left(\int_{\mD_M}|\na\bl{w}|^2dx+\int_{\Sigma_L}\left|\bl{w}(x_h,-Z/2)\right|^2dS+\int_{\Sigma_R}\left|\bl{w}(x_h,Z/2)\right|^2dS\right).
\ee
Meanwhile, using the trace theorem in $\Sigma_L\times[-Z,-Z/2]$ and Lemma \ref{P2}, one derives
\be\label{PCC2}
\begin{split}
\int_{\Sigma_L}\left|\bl{w}(x_h,-Z/2)\right|^2dS&\leq C_\mD\left(\int_{\Sigma_L\times[-Z,-Z/2]}\left|\bl{w}\right|^2dx+\int_{\Sigma_L\times[-Z,-Z/2]}\left|\na\bl{w}\right|^2dx\right)\\
&\leq C_\mD\int_{\Sigma_L\times[-Z,-Z/2]}\left|\na\bl{w}\right|^2dx.
\end{split}
\ee
Similarly, one derives that
\be\label{PCC3}
\int_{\Sigma_R}\left|\bl{w}(x_h,Z/2)\right|^2dS\leq C_\mD\int_{\Sigma_R\times[Z/2,Z]}\left|\na\bl{w}\right|^2dx.
\ee
Substituting \eqref{PCC2}--\eqref{PCC3} in the right hand side of \eqref{PCC1}, one deduces
\[
\int_{\mD_M}|\bl{w}|^2dx\leq C_\mD\int_{D_\zeta}|\na\bl{w}|^2dx.
\]
This finishes the estimate in $\mD_M$. Noting that the rest parts of the domain $D_\zeta$ are a union of two straight truncated pipes, estimates in $\mD_\zeta-\mD_M$ are direct conclusions of Lemma \ref{P2}. This concludes the validity of \eqref{TORPIPEPOIN}.

\qed

The following asymptotic estimate of a function that satisfies an ordinary differential inequality will be useful in our further proof. To the best of the authors' knowledge, it was originally derived by Ladyzhenskaya-Solonnikov in \cite{Lady-Sol1980}.
\begin{lemma}\label{LEM2.3}
Let $Y(\zeta)\nequiv 0$ be a nondecreasing nonnegative differentiable function satisfying
\be\label{cond}
Y(\zeta)\leq\Psi(Y'(\zeta)),\q\forall\zeta>0.
\ee
Here $\Psi:\,[0,\infty)\to[0,\infty)$ is a monotonically increasing function with $\Psi(0)=0$ and there exists $C,\,\tau_1>0$, $m>1$, such that
\be\label{cond1}
\Psi(\tau)\leq C\tau^m,\q\forall\tau>\tau_1.
\ee
Then
\be\label{LEM2.333}
\liminf_{\zeta\to+\infty}\zeta^{-\f{m}{m-1}}Y(\zeta)>0.
\ee
\end{lemma}

\pf Since $Y$ is not identically zero, there exists $\zeta_0>0$ such that $Y(\zeta_0):=Y_0>0$. Using the monotonicity of $\Psi$, one knows that
\[
Y'(\zeta_0)\geq\Psi^{-1}(Y_0):=\eta_0>0.
\]
And therefore we have
\[
Y(\zeta)\geq Y_0+\eta_0(\zeta-\zeta_0)\q \text{for } \zeta\geq \zeta_0.
\]
From \eqref{cond}, we see that
\[
Y'(\zeta)\geq\Psi^{-1}(Y(\zeta))\geq\Psi^{-1}(Y_0+\eta_0(\zeta-\zeta_0)) \q \text{for } \zeta\geq \zeta_0.
\]
When $\zeta\to +\i$, we can deduce that $\Psi^{-1}(Y_0+\eta_0(\zeta-\zeta_0))\to+\i$, otherwise if $\Psi^{-1}(Y_0+\eta_0(\zeta-\zeta_0))\to A<+\i$,
\bes
Y_0+\eta_0(\zeta-\zeta_0)=\Psi\lt(\Psi^{-1}(Y_0+\eta_0(\zeta-\zeta_0))\rt)\leq \Psi(A)<+\i,
\ees
which is invalid as $\zeta\to +\i$.

So we see that there exists $\zeta_1\geq\zeta_0$ such that $Y'(\zeta)\geq\tau_1$ for any $\zeta\geq\zeta_1$. Then from \eqref{cond} and \eqref{cond1}, we have for $\zeta\geq\zeta_1$,
\[
Y(\zeta)\leq C(Y'(\zeta))^m.
\]
Integrating the above inequality on $[\zeta_1,\infty)$, one concludes \eqref{LEM2.333}.

\qed

At the end of this section, we introduce the following results which focus on the $W^{1,3}$-weak solution and the $H^m$-strong solution of the linear Stokes equations on \emph{bounded domains} with the Navier-slip boundary condition.

\begin{lemma}[See \cite{AACG:2021JDE}, Corollary 5.7]\label{LEM52}
Let $\Omega$ be a bounded smooth domain, $\mathbf{f} \in L^{\f{3}{2}}(\Omega), {\mathbf{F}} \in L^{3}(\Omega)$, and ${\mathbf{h}} \in {W}^{-\frac{1}{3}, 3}(\p\O)$. Then, the Stokes problem
\[
\left\{
\begin{array}{*{2}{ll}}
-\Dl \mathbf{v}+\na P=\mathrm{div}\,\mathbf{F}+\mathbf{f}, \q\na\cdot \mathbf{v}=0,\q\q&\text{in }\O;\\[2mm]
2(\mathbb{S}\mathbf{v}\cdot\bl{n})_{\mathrm{tan}}+\al \mathbf{v}_{\mathrm{tan}}=\bl{h},\q \mathbf{v}\cdot\bl{n}=0, \q\q&\text{on }\p\O,
\end{array}
\right.
\]
has a unique solution $(\mathbf{v}, P) \in W^{1, 3}(\Omega) \times L^{3}(\Omega)$ which satisfies the estimate:
\be\label{ESTLEM27}
\|\mathbf{v}\|_{{W}^{1, 3}(\Omega)}+\|P\|_{L^{3}(\Omega)} \leq C_{\alpha, \Omega }\left(\|\mathbf{f}\|_{{L}^{\f{3}{2}}(\Omega)}+\|\mathbf{F}\|_{{L}^{3}(\Omega)}+\|\mathbf{h}\|_{{W}^{-\frac{1}{3}, 3}(\p\O)}\right).
\ee
\end{lemma}

\qed

\begin{lemma}[See \cite{AACG:2021JDE}, Theorem 4.5, and \cite{Ghosh:2018PHD}, Theorem 2.5.10]\label{LEM59}
Let $\Omega$ be a bounded smooth domain, $m\in\mathbb{N}$, $\mathbf{f}\in H^{m}(\Omega)$, and $\mathbf{h} \in {H}^{m+\frac12}(\p\O)$. Then, the solution of the Stokes problem
\[
\left\{
\begin{array}{*{2}{ll}}
-\Dl \mathbf{v}+\na P=\mathbf{f}, \q\na\cdot \mathbf{v}=0,\q\q&\text{in }\O;\\[2mm]
2(\mathbb{S}\mathbf{v}\cdot\bl{n})_{tan}+\al \mathbf{v}_{tan}=\bl{h},\q \mathbf{v}\cdot\bl{n}=0, \q\q&\text{on }\p\O,
\end{array}
\right.
\]
satisfies $(\mathbf{v}, P) \in {H}^{m+2}(\Omega) \times H^{m+1}(\Omega)$. Also, it enjoys the following estimate
\be\label{ESTLEM28}
\|\mathbf{v}\|_{{H}^{m+2}(\Omega)}+\|P\|_{H^{m+1}(\Omega)}\leq C_{\al,\Omega}\left(\|\mathbf{f}\|_{{H}^{m}(\Omega)}+\|\mathbf{h}\|_{{H}^{m+\frac{1}{2}}(\p\O)}\right).
\ee
\end{lemma}

\qed

\section{Existence}\label{SECE}
\subsection{On Poiseuille flows in pipes $\bl{\mD_L}$ and $\bl{\mD_R}$}
\q In this subsection,  we introduce Poiseuille flows in pipes ${\mD_L}$ and ${\mD_R}$, which are solutions of system \eqref{GL1} in $\Sigma_i\times\bR$ ($i=L$ or $R$).  We will drop the index $i$ for convenience in this subsection. To find a Poiseuille flow $\bl{g}_\Phi$ in $\Sigma\times\bR$ with a given flux $\Phi$, one needs to find a function $g_\Phi:\Sigma\to\mathbb{R}$, such that
\be\label{HP1}
\left\{
\begin{aligned}
&\bl{g}_\Phi=g_\Phi\bl{e_3}, \\
&-\Dl_hg_\Phi(x_h)=Constant,\q\q\q\text{in }\Sigma,\\
&\f{\p g_\Phi}{\p\bar{\bl{n}}}=-\al g_\Phi,\hskip 1.9cm\q\q\q\text{on }\p\Sigma,\\
&\int_{\Sigma}g_\Phi(x_h)dx_h=\Phi. \\
\end{aligned}
\right.
\ee
Here and below, we assume $\Phi\geq0$ with out loss of generality.
\begin{figure}[H]
\centering
\includegraphics[scale=0.53]{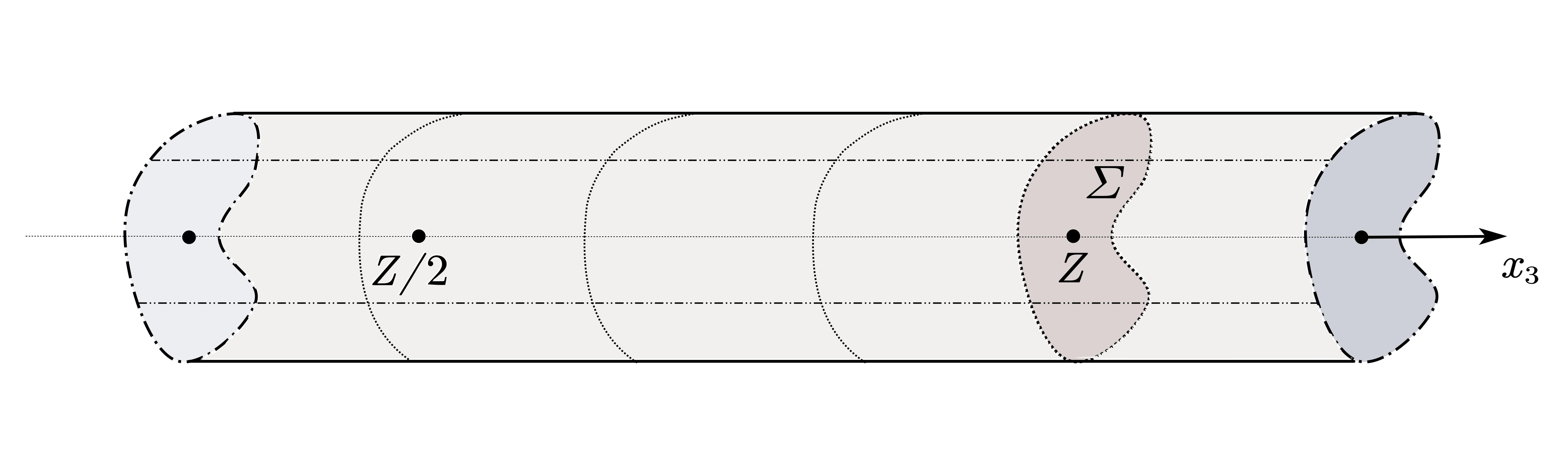}
\caption{Straight infinite pipe $\Sigma\times\mathbb{R}$.}
\end{figure}
\begin{remark}
If $\Sigma$ is the unit disk in $\mathbb{R}^2$, one has the following exact formula of $\bl{g}_\Phi$ :
\[
\bl{g}_\Phi(x)=\f{2(\al+2)\Phi}{(\al+4)\pi}\lt(1-\f{\al}{\al+2}|x_h|^2\rt)\boldsymbol{e_3},\q\text{with its pressure } p_\Phi(x)=-\f{8\al\Phi}{(\al+4)\pi} x_3,
\]
which could be considered as a generalization of the Hagen-Poiseuille flow under the no-slip boundary condition ($\al\to +\i$).
\end{remark}

\qed

The existence and uniqueness of $g_\Phi$ in \eqref{HP1} could be derived by routine methods of elliptic equations in a bounded smooth domain with the Robin boundary condition, here we omit the details. We only derive an $H^m$ estimate of $g_\Phi$ in terms of $\Phi$ by the scaling technique.

 By the linearity of the problem \eqref{HP1}, one considers the classical Poisson equation subject to the Robin boundary condition:

\be\label{ROB}
\left\{
\begin{aligned}
-\Dl_h\varphi(x_h)&=1,\q\text{in }\Sigma,\\
\f{\p\varphi}{\p\bar{\bl{n}}}+\al\varphi&=0,\q\text{on }\p\Sigma.\\
\end{aligned}
\right.
\ee
The existence and uniqueness of the problem \eqref{ROB} are classical. We refer readers to \cite{GT1998} for details. Moreover, the following estimate of $\varphi$ can be derived by classical results:
\[
\|\varphi\|_{H^{m}(\Sigma)}\leq C_{\al,m,\Sigma}.
\]
Multiplying \eqref{ROB}$_1$ by $\varphi$ and integration by parts, we also note that:
\be\label{POS}
\int_{\Sigma}\varphi dx_h=-\int_{\Sigma}\varphi\Dl_h\varphi dx=\int_{\Sigma}|\nabla_h\varphi|^2dx+\int_{\p\Sigma}\al|\varphi|^2dS_h:=C_P>0.
\ee
Thus one concludes that
\[
g_\Phi(x_h)=\f{\Phi}{C_P}\varphi(x_h)
\]
solves the problem \eqref{HP1}. Then
\be\label{EHP}
\bali
\|g_\Phi\|_{H^{m}(\Sigma)}=\f{\Phi}{C_P}\|\varphi(x_h)\|_{H^{m}(\Sigma)}\leq C_{\al,m,\Sigma}\Phi,\q\forall m\in\mathbb{N}
\eali
\ee
where $C_{\al,m,\Sigma}>0$ is a constant independent of $\Phi$. Later, for $i=L$ or $R$, we denote $\bl{g}^i_{\Phi}=g^i_{\Phi}\bl{e_3}$ the Poiseuille flows in $\Sigma_i\times\bR$.

\subsection{Construction of the profile vector}\label{SEC3.2}

\q\ In this subsection, we are devoted to the construction of a smooth divergence-free vector $\bl{a}$, which satisfies the Navier-slip boundary condition \eqref{NBC}. Meanwhile, the vector $\bl{a}$ equals to $\bl{g}^L_{\Phi}$ in far left of $\mD$, and it is identical with $\bl{g}^R_{\Phi}$ in far right of $\mD$. Here is the result:

\begin{proposition}\label{Prop}
There exists a smooth vector field $\bl{a}(x)$ which enjoys the following properties
\begin{itemize}
\setlength{\itemsep}{-2 pt}
\item[(i).] $\bl{a}\in C^\i(\overline{\mD})$, and $\na\cdot \bl{a}=0$ in $\mD$;

\item[(ii).] $2(\mathbb{S}\bl{a}\cdot\bl{n})_{\mathrm{tan}}+\al \bl{a}_{\mathrm{tan}}=0$, and $\bl{a}\cdot\bl{n}=0$ on $\p\mD$;

\item[(iii).] $\bl{a}=\bl{g}^L_{\Phi}$ in $\mD\cap\{x\in\mathbb{R}^3:\,x_3\leq -Z\}$, and $\bl{a}=\bl{g}^R_{\Phi}$ in $\mD\cap\{x\in\mathbb{R}^3:\,x_3\geq Z\}$;

\item[(iv).] $\|\bl{a}\|_{H^m(\mD_M)}\leq C_{\al,m,\mD}\Phi$,\q for any $m\in\mathbb{N}$.
\end{itemize}
\end{proposition}

\qed

\pf Recalling the assumption of the domain $\mD$, one notices that there exists a smooth domain $\Sigma'\subset\mathbb{R}^2$ such that $\Sigma'\times\mathbb{R}\subset\subset\mD$ (\emph{which means $\Sigma'\times\mathbb{R}\subset\mD$ and $\text{dist}(\Sigma'\times\mathbb{R},\p\mD)\geq\ve_0>0$}). Let $h=h(x_h)$ be a smooth function supported on $\Sigma'$ which satisfies
\[
\int_{\Sigma'}h(x_h)dx_h=\Phi.
\]
By a scaling, we can assume that
\be\label{hphi}
\|h\|_{H^m(\Sigma')}\leq C_m\Phi,\q\forall m\in\mathbb{N}.
\ee

Let $\eta=\eta(x_3)$ be smooth cut-off functions such that
\be\label{ETA}
\eta(x_3)=\left\{
\begin{aligned}
&1,\q\text{for}\q x_3>Z;\\
&0,\q\text{for}\q x_3<Z/2.\\
\end{aligned}
\right.
\ee
Now we define the vector $\bl{a}$ by

\be\label{CONSa}
\bl{a}=\lt(
\begin{matrix}
A^R_1(x_h)\eta'(x_3)-A^L_1(x_h)\eta'(-x_3) \\
A^R_2(x_h)\eta'(x_3)-A^L_2(x_h)\eta'(-x_3)\\
h(x_h)+(g^R_{\Phi}-h(x_h))\eta(x_3)+(g^L_{\Phi}-h(x_h))\eta(-x_3)
\end{matrix}
\rt)^T,
\ee
where the 2D vector function $\bl{A}^i:=(A^i_1,A^i_2)$ ($i=L$ or $R$)  solves the following partial differential equation in $\Sigma_i$
\be\label{Aeq}
\text{div}_h \bl{A}^i(x_h)=h(x_h)-g^i_{\Phi}(x_h),\q\text{in}\q \Sigma_i,
\ee
subject to the two dimensional Navier-slip boundary condition

\be\label{Abc}
\left\{
\begin{aligned}
&2(\mathbb{S} \bl{A}^i\cdot\bar{\bl{n}})_{\mathrm{tan}}+\al \bl{A}^i_{\mathrm{tan}}=0,\\
&\bl{A}^i\cdot\bar{\bl{n}}=0,\\
\end{aligned}
\right.\q\text{on}\q\p\Sigma_i,
\ee
where $\bar{\bl{n}}=(n_1(x_h), n_2(x_h))$ is the unit outward normal vector on $\p\Sigma_i$.


Now let us verify the validity of the above construction.

First, combining \eqref{CONSa} and \eqref{Aeq}, direct computation shows the divergence-free property of $\bl{a}$. The smoothness of $\bl{a}$ follows from the smoothness of $h$, $\eta$, $g_\Phi^L$ and $g_\Phi^R$ which are provided in their definitions, together with the smoothness of $\bl{A}^L$ and $\bl{A}^R$ which will be derived below.

Second, concerning the validity of boundary condition (ii) in Proposition \ref{Prop}, we first see that
\[
\bl{a}=\left\{
\begin{array}{*{2}{ll}}
g^L_\Phi(x_h)\bl{e_3},&\q \text{ in }\q \mD\cap\{x:\,x_3\leq -Z\};\\[2mm]
g^R_\Phi(x_h)\bl{e_3},&\q \text{ in }\q \mD\cap\{x:\,x_3\geq Z\};\\[2mm]
h(x_h)\bl{e_3},&\q \text{ in }\q \mD\cap\{x:\,|x_3|\leq Z/2\}.
\end{array}
\right.
\]
Due the fact that $g^i_\Phi\bl{e_3}$ ($i=L,R$) is the Poiseuille flow in \eqref{HP1} which satisfies the same Navier-slip boundary condition, and the auxiliary function $h(x_h)$ is compactly supported in each \emph{cross section} of $\mD$, we see $\bl{a}$ satisfies the Navier-slip boundary condition on $\p\mD\cap\{x:\,|x_3|\leq Z/2\ \text{or}\ |x_3|\geq Z\}$.  For the remaining part $\p\mD\cap\{x:\,Z/2\leq|x_3|\leq Z\}$, the unit outer normal vector enjoys the following from
 \[
 \bl{n}=(\bar{\bl{n}},0)=(n_1(x_h), n_2(x_h),0),\q\text{on } \p\mD\cap\{x:\,Z/2\leq|x_3|\leq Z\},
 \]
which is independent with $x_3$ variable. Recalling \eqref{NBCM1}, the Navier-slip boundary condition
\[
\left\{
\begin{aligned}
&2(\mathbb{S} \bl{a}\cdot{\bl{n}})_{\mathrm{tan}}+\al \bl{a}_{\mathrm{tan}}=0,\\
&\bl{a}\cdot{\bl{n}}=0,\\
\end{aligned}
\right.\q\text{on}\q\p\mD\cap\{x:\,Z/2\leq|x_3|\leq Z\}
\]
enjoys the following form  in the orthogonal curvilinear coordinates on the boundary:
\be\label{BCBCC}
\left\{
\begin{aligned}
&\p_{\boldsymbol{n}}a_{\tau_1}=\left(\kappa_1(x)-\al\right)a_{\tau_1},\\
&\p_{\boldsymbol{n}}a_3=-\al a_3,\\
&a_n=0.\\
\end{aligned}
\right.\q\q\q\text{on}\q\p\mD\cap\{x:\,Z/2\leq|x_3|\leq Z\}.
\ee
Therefore, noting that the cut-off function $\eta$ depends only on $x_3$ variable, \eqref{BCBCC}$_{1}$ and \eqref{BCBCC}$_{3}$ are guaranteed by \eqref{Abc}$_1$ and \eqref{Abc}$_2$, respectively. Moreover, by direct calculation,
\[
\begin{split}
\p_{\bl{n}}a_3=\eta(x_3)\p_{\bl{n}}g_{\Phi}^R+\eta(-x_3)\p_{\bl{n}}g_{\Phi}^L=-\al\eta(x_3)g_{\Phi}^R-\al\eta(-x_3)g_{\Phi}^L=&-\al a_3,\\
&\text{on}\q\p\mD\cap\{x:\,Z/2\leq|x_3|\leq Z\},
\end{split}
\]
one proves \eqref{BCBCC}$_2$. This concludes the validity of item (ii) in Proposition \ref{Prop}.

Third, property (iii) in this proposition follows directly by the definition of $\bl{a}$ in \eqref{CONSa}.

Finally, we derive the $H^m$ estimate on $\bl{a}$ in $\mD_M$. Using \eqref{CONSa}, we see
\be\label{EOFA627}
\bali
\|\bl{a}\|_{H^m(\mD_M)}&\ls \sum_{i=L,R}\left\|\bl{A}^i\right\|_{H^m(\Sigma_i)}+\|h\|_{H^m(\Sigma')}+ \sum_{i=L,R}\|{g}^i_{\Phi}\|_{H^m(\Sigma_i)}\\
 &\ls_{\al,m,\mD} \Phi+ \sum_{i=L,R}\left\|\bl{A}^i\right\|_{H^m(\Sigma_i)},
\eali
\ee
where the last inequality follows from the estimates in \eqref{EHP} and \eqref{hphi}. Now it remains to show the $H^m$ estimate of 2D vectors $\bl{A}^i$, $i=L,\,R$, by solving the boundary value problem \eqref{Aeq}--\eqref{Abc}, which is derived in the following lemma:
\begin{lemma}
Problem \eqref{Aeq}--\eqref{Abc} has a smooth solution $\bl{A}^i\in C^\i(\overline{\Sigma})$ satisfying
\be\label{ESTLEM3.3}
\left\|\bl{A}^i\right\|_{H^m(\Sigma_i)}\leq C_{\al,m,\Sigma_i}\Phi,\q\forall m\in\mathbb{N}.
\ee
\end{lemma}

\pf For notation of simplicity, we will omit the index $L$ or $R$ in the following proof if no ambiguity is caused. Using the Helmholtz-Weyl decomposition, we can split $\bl{A}$ into
\be\label{EEEE}
\bl{A}:=\nabla_h \phi+\bl{G}.
\ee
Here $\phi=\phi(x_h)$ is a scalar function which satisfies
\be\label{A1}
\begin{cases}
\Dl_h\phi=h(x_h)-g_{\Phi}(x_h),&\q\text{in}\q\Sigma;\\[2mm]
\f{\p\phi}{\p\bar{\bl{n}}}=0,&\q\text{on}\q\p\Sigma;\\[2mm]
\int_\Sigma\phi dx_h=0.
\end{cases}
\ee
By the definition of the auxiliary function $h(x_h)$, one has $h-g_\Phi$  satisfies the following compatibility condition:
\[
\int_{\Sigma}\big(h(x_h)-g_\Phi(x_h)\big)dx_h=0.
\]
Thus classical theory of Poisson equations indicates the solvability and regularity $\phi\in C^\i(\overline{\Sigma})$  of problem \eqref{A1}. And $\phi$ satisfies
\be\label{phic}
\|\phi\|_{H^{m+2}(\Sigma)}\leq C_{m,\Sigma}\|h-g_{\Phi}\|_{H^m(\Sigma)}\leq C_{m,\Sigma} \Phi,\q\forall m\in\mathbb{N}.
\ee

It remains to construct the smooth vector $\bl{G}$ in \eqref{EEEE}. Notice that $\bl{G}$ should satisfy

\be\label{AC}
\left\{
\begin{array}{*{2}{ll}}
\text{div}_h\bl{G}=0,&\q\text{in}\q\Sigma,\\[2mm]
2(\mathbb{S} \bl{G}\cdot\bar{\bl{n}})_{\mathrm{tan}}+\al \bl{G}_{\mathrm{tan}}=2(\mathbb{S} (\na\phi)\cdot\bar{\bl{n}})_{\mathrm{tan}}+\al (\na\phi)_{\mathrm{tan}},&\q\text{on}\q\p\Sigma,\\[2mm]
\bl{G}\cdot\bar{\bl{n}}=0,&\q\text{on}\q\p\Sigma.
\end{array}
\right.
\ee
There is too much space for us to construct a solution $G$ satisfying \eqref{AC} such that $\|G\|_{H^m(\Sigma)}\leq C_{\al,m,\Sigma}\Phi$.  For example, we can choose $(G, \pi)$ to be the pair of solution to the following linear Stokes equations with the Navier-slip boundary condition:
\[
\left\{
\begin{array}{*{2}{ll}}
-\Dl_h \bl{G}+\na\pi=0,\q \text{div}_h\bl{G}=0,&\q\text{in}\q\Sigma,\\[2mm]
2(\mathbb{S} \bl{G}\cdot\bar{\bl{n}})_{\mathrm{tan}}+\al \bl{G}_{\mathrm{tan}}=2(\mathbb{S} (\na\phi)\cdot\bar{\bl{n}})_{\mathrm{tan}}+\al (\na\phi)_{\mathrm{tan}},&\q\text{on}\q\p\Sigma,\\[2mm]
\bl{G}\cdot\bar{\bl{n}}=0,&\q\text{on}\q\p\Sigma.
\end{array}
\right.
\]

From \cite[Theorem 4.5]{AACG:2021JDE} or \cite[Theorem 2.5.10]{Ghosh:2018PHD}\footnote{Strictly speaking, theorems in \cite{AACG:2021JDE,Ghosh:2018PHD} are derived for 3D linear Stokes systems. However, their methods are also valid for related 2D problems. See the introduction part of \cite{Ghosh:2018PHD}.}, we have the following estimate of $\bl{G}$

\be\label{phic1}
\bali
\|\bl{G}\|_{H^{m+2}(\Sigma)}\leq& C_{\al,m,\Sigma}\|2(\mathbb{S} (\na\phi)\cdot\bar{\bl{n}})_{\mathrm{tan}}+\al (\na\phi)_{\mathrm{tan}}\|_{H^{m+1/2}(\p\Sigma)}\\
 \leq& C_{\al,m,\Sigma}\|\phi\|_{H^{m+3}(\Sigma)}\leq C_{\al,m,\Sigma}\Phi,
\eali
\ee
where at the last line of the above inequality, we have used the trace theorem and \eqref{phic}. Then \eqref{ESTLEM3.3} is proved by combining \eqref{phic} and \eqref{phic1}.

\qed

\begin{remark}
Combining estimates \eqref{EHP}, \eqref{EOFA627} and \eqref{phic1} above, the following global $W^{1,\i}$-estimate of $\bl{a}$ is a direct conclusion of the Sobolev imbedding:
\be\label{EOFAA6277}
\|\bl{a}\|_{W^{1,\i}(\mD)}\leq C_{\al,\mD}\Phi.
\ee
\end{remark}

\qed

\subsection{Proof of the existence}\label{SEC3627}
\q\ In this section we study the solvability of the generalized Leray's problem subject to the Navier-slip boundary condition. Consider the asymptotic behavior of the prescribed weak solution in Definition \ref{weaksd1}, we write
\be\label{EEPPSS}
\bl{u}=\bl{v}+\bl{a},
\ee
where $\bl{a}$ is constructed in the previous subsection. Therefore, the generalized Leray's problem \eqref{GL1}--\eqref{GL3} has the following equivalent form in the viewpoint of $\bl{v}$:
\begin{problem}[Modified problem]\label{PP2}
Find $(\bl{v},p)$ such that
\be\label{EQU}
\left\{
\begin{aligned}
&\bl{v}\cdot\nabla \bl{v}+\bl{a}\cdot\nabla \bl{v}+\bl{v}\cdot\nabla \bl{a}+\na p-\Dl \bl{v}=\Dl \bl{a}-\bl{a}\cdot\nabla \bl{a},\\
&\na\cdot \bl{v}=0,\\
\end{aligned}
\right.\q\text{in }\mD,
\ee
subject to the Navier-slip boundary condition
\be\label{BDR}
\left\{
\begin{aligned}
&2(\mathbb{S}\bl{v}\cdot\bl{n})_{\mathrm{tan}}+\al \bl{v}_{\mathrm{tan}}=0,\\
&\bl{v}\cdot\bl{n}=0,\\
\end{aligned}
\right.\q\text{on }\p\mD,
\ee
with the asymptotic behavior as $|x_3|\to\i$
\be\label{ASY}
\bl{v}(x)\to \bl{0},\q\text{as}\q |x_3|\to \i.
\ee
\end{problem}

\qed

Substituting the expression \eqref{EEPPSS} in the weak formulation \eqref{weaksd}, we arrive at the following weak formulation of $\bl{v}$:
\begin{definition}\label{defws}
Let $\bl{a}$ be a smooth vector satisfying the properties stated in Proposition \ref{Prop}. We say that $\bl{v}\in \mathcal{H}_\sigma(\mD)$ is a weak solution of Problem \ref{PP2} if
\be\label{EQU111}
\begin{split}
&2\int_{\mD}\mathbb{S}\bl{v}:\mathbb{S}\bl{\varphi} dx+\al\int_{\p\mD}\bl{v}_{\mathrm{tan}}\cdot\bl{\varphi}_{\mathrm{tan}}dS+\int_{\mD}\bl{v}\cdot\nabla \bl{v} \cdot \bl{\varphi} dx+\int_{\mD}\bl{v}\cdot\nabla \bl{a}\cdot\bl{\varphi} dx\\
+&\int_{\mD}\bl{a}\cdot\nabla \bl{v}\cdot\bl{\varphi} dx=\int_{\mD}\big(\Dl \bl{a}-\bl{a}\cdot\nabla \bl{a}\big)\cdot\bl{\varphi} dx
\end{split}
\ee
holds for any vector-valued function $\bl{\varphi}\in \mathcal{H}_\sigma(\mD)$.
\end{definition}

\qed

To establish the existence of the weak solution defined in Definition \ref{defws}, we first introduce the following Brouwer's fixed point theorem. It could be found in \cite{Lions1969}. See also \cite[Lemma IX.3.1]{Galdi2011}.
\begin{lemma}\label{FUNC}
Let $P$ be a continuous operator which maps $\mathbb{R}^N$ into itself, such that for some $\rho>0$
\[
P({\xi})\cdot{\xi}\geq 0 \quad \text { for all } {\xi} \in\mathbb{R}^n \text { with }|{\xi}|=\rho.
\]
Then there exists ${\xi}_{0} \in\mathbb{R}^N$ with $|{\xi}_{0}| \leq \rho$ such that ${P}({\xi}_{0})=0$.
\end{lemma}

\qed

Now, we go to the existence theorem.
\begin{theorem}\label{THM3.8}
There is a constant $\Phi_0>0$ depending on $\al$ and the curvature of $\p\mD$ such that if $\Phi\leq\Phi_0$, then Problem \ref{PP2} admits at least one weak solution $(\bl{v},p)\in \mathcal{H}_\sigma(\mD)\times L^2_{\mathrm{loc}}(\overline{\mD})$, with
\be\label{EOFV627}
\|\bl{v}\|_{H^1(\mD)}\leq C_{\al,\mD}\Phi.
\ee
\end{theorem}

\qed

\begin{remark}
The weak solution satisfies a generalized version of \eqref{ASY}. Actually,
it follows from the trace inequality (\cite[Theorem II.4.1]{Galdi2011}) that
\[\int_{\Sigma_R}|\bl{v}(x_h,x_3)|^2dx_h\leq C\int_{z>x_3}\int_{\Sigma_R}(|\bl{v}|^2+|\nabla \bl{v}|^2)(x_h,z)dx_hdz\,,\]
where the constant $C$ is independent of $x_3$ variable. This implies
\[
\int_{\Sigma_R}|\bl{v}(x_h,x_3)|^2dx_h\to0,\quad\text{as $x_3\to+\infty$}\,.
\]
The case when $x_3\to-\i$ is similar.
\end{remark}

\qed

Now we are ready to provide the proof of Theorem \ref{THM3.8}.

\subsubsection{Constructing the velocity field by the Galerkin method}
\q\ Using the Galerkin method, we first construct an approximate solution and then pass to the limit by compactness arguments. Recalling
\[
\bl{X}:=C^\infty_{\sigma,c}(\overline{\mD};\,\mathbb{R}^3)=\big\{\bl{\varphi}\in C^\infty_c(\overline{\mD}\,;\,\mathbb{R}^3):\, \na\cdot \bl{\varphi}=0,\ \bl{\varphi}\cdot\bl{n}\big|_{\p\mD}=0\big\},
\]
and $\{\bl{\varphi}_k\}_{k=1}^\i\subset\bl{X}$ be an unit orthonormal basis of $\mathcal{H}_\sigma(\mD)$, that is:
\[
\langle \bl{\varphi}_i,\bl{\varphi}_j\rangle_{H^1(\mD)}=\begin{cases}
1,&\q\text{if }\q i=j;\\
0,&\q\text{if }\q i\neq j,\\
\end{cases}
\]
$\forall i,j\in\mathbb{N}$. Now we construct an approximation of $\bl{v}$ of the form
\[
\bl{v}_N(x)=\sum_{i=1}^Nc_i^N\bl{\varphi}_i(x).
\]
To determine $\bl{v}_N$, one tests the weak formulation \eqref{EQU111} by $\bl{\varphi}_i$, with $i=1,2,...,N$. This indicates that
\[
\begin{split}
&2\sum_{i=1}^Nc_i^N\int_{\mD}\mathbb{S}\bl{\varphi}_i:\mathbb{S}\bl{\varphi}_j dx+\al\sum_{i=1}^N c_i^N\int_{\p\mD} (\bl{\varphi}_i)_{\mathrm{tan}}(\bl{\varphi}_j)_{\mathrm{tan}}dS+\sum_{i,k=1}^Nc_i^Nc_k^N\int_{\mD}\bl{\varphi}_i\cdot\nabla \bl{\varphi}_k\cdot \bl{\varphi}_jdx\\
+&\sum_{i=1}^N\int_{\mD}\bl{\varphi}_i\cdot\nabla \bl{a}\cdot \bl{\varphi}_jdx+\sum_{i=1}^Nc_i^N\int_{\mD}\bl{a}\cdot\nabla \bl{\varphi}_i\cdot \bl{\varphi}_jdx=\int_{\mD}\big(\Dl \bl{a}-\bl{a}\cdot\nabla \bl{a}\big)\cdot \bl{\varphi}_jdx,\q\forall j=1,2,...,N.
\end{split}
\]
As we see, this is a system of nonlinear algebraic equations of $N$-dimensional vector
\[
\bl{c}^N:=(c_1^N,c_2^N,...,c_N^N).
\]

We shall solve the above system by Lemma \ref{FUNC} (Brouwer's fixed point theorem). To this end, we denote ${P}:\,\mathbb{R}^N\to\mathbb{R}^N$ such that
\bes
\begin{split}
\big({P}(\bl{c}^N)\big)_j=&2\sum_{i=1}^Nc_i^N\int_{\mD}\mathbb{S}\bl{\varphi}_i:\mathbb{S}\bl{\varphi}_j dx+\al\sum_{i=1}^N c_i^N\int_{\p\mD} (\bl{\varphi}_i)_{\mathrm{tan}}\cdot(\bl{\varphi}_j)_{\mathrm{tan}}dS+\sum_{i,k=1}^Nc_i^Nc_k^N\int_{\mD}\bl{\varphi}_i\cdot\nabla \bl{\varphi}_k\cdot \bl{\varphi}_jdx\\
&+\sum_{i=1}^N\int_{\mD}\bl{\varphi}_i\cdot\nabla \bl{a}\cdot \bl{\varphi}_jdx+\sum_{i=1}^Nc_i^N\int_{\mD}\bl{a}\cdot\nabla \bl{\varphi}_i\cdot \bl{\varphi}_jdx-\int_{\mD}\big(\Dl \bl{a}-\bl{a}\cdot\nabla \bl{a}\big)\cdot \bl{\varphi}_jdx,\\
&\hskip 10cm\q\forall j=1,2,...,N.
\end{split}
\ees
Clearly, one observes that ${P}$ is continuous. Then we check that
\[
{P}(\bl{c}^N)\cdot\bl{c}^N=\un{2\int_{\mD}|\mS \bl{v}_N|^2dx+\al\int_{\p\mD}|(\bl{v}_N)_{\mathrm{tan}}|^2dS}_{I_1}+\un{\int_{\mD}\left((\bl{v}_N+\bl{a})\cdot\nabla(\bl{v}_N+\bl{a})\right)\cdot \bl{v}_Ndx}_{I_2}-\un{\int_{\mD}\bl{v}_N\cdot\Dl \bl{a}dx}_{I_3}.
\]
First, we estimate the term $I_1$. We will show that there exists a constant $C_{\al,\mD}$, depending on $\al$ and $\mD$, such that
\be\label{termi1}
I_1\geq C_{\al,\mD}\int_{\mD}|\nabla \bl{v}_N|^2dx.
\ee

By the definition of the stress tensor and integration by parts, one notices that
\be\label{EM111}
\begin{split}
\int_{\mD}|\mS \bl{v}_N|^2dx&=\f{1}{2}\int_{\mD}|\na \bl{v}_N|^2dx+\f{1}{2}\sum_{i,j=1}^3\int_{\mD}\p_{x_i}(\bl{v}_N)_j\p_{x_j}(\bl{v}_N)_idx\\
&=\f{1}{2}\int_{\mD}|\na \bl{v}_N|^2dx+\f{1}{2}\un{\sum_{i,j=1}^3\int_{\p\mD}(\bl{v}_N)_j\p_{x_j}(\bl{v}_N)_in_idS}_{I_{11}}-\f{1}{2}\un{\int_{\mD}\bl{v}_N\cdot\nabla\text{div }(\bl{v}_N)dx}_{I_{12}}.
\end{split}
\ee
Here the term $I_{12}$ vanishes due to $\bl{v}_N$ is divergence-free. Noting that $\bl{v}_N\cdot\bl{n}\equiv0$ on $\p\mD$, it implies
\[
I_{11}=\int_{\p\mD}\bl{v}_N\cdot\big(\nabla\left(\bl{v}_N\cdot\bl{n}\right)-\bl{v}_N\cdot\nabla\bl{n}\big)dS=-\sum_{i,j=1}^3\int_{\p\mD}(\bl{v}_N)_j\p_{x_j}n_i(\bl{v}_N)_jdS.
\]
Thus $I_{11}$ can be bounded by
\[
|I_{11}|\leq C_{\mD}\int_{\p\mD}|(\bl{v}_N)_{\mathrm{tan}}|^2dS,
\]
where $C_\mD>0$ is a universal constant, depending only on $\p\mD$. Inserting the above calculations for $I_{11}$ and $I_{12}$ in \eqref{EM111}, one arrives at
\[
I_1\geq2\int_{\mD}|\mS \bl{v}_N|^2dx\geq\int_{\mD}|\na \bl{v}_N|^2dx-C_{\kappa}\int_{\p\mD}|(\bl{v}_N)_{\mathrm{tan}}|^2dS,
\]
that is
\[
\frac{\alpha}{C_{\mD}}I_1\geq \frac{\alpha}{C_{\mD}}\int_{\mD}|\na \bl{v}_N|^2dx-\alpha\int_{\p\mD}|(\bl{v}_N)_{\mathrm{tan}}|^2dS\,.
\]
Hence, we deduce that
\[
I_1+\frac{\alpha}{C_{\mD}}I_1\geq \frac{\alpha}{C_{\mD}}\int_{\mD}|\na \bl{v}_N|^2dx\,.
\]
This indicates \eqref{termi1}.

Next, we turn to the estimate of $I_2$. Using integration by parts, together with the divergence-free property of $\bl{v}_N$ and $\bl{a}$, one knows that
\[
I_2=\un{\int_{\mD}\bl{v}_N\cdot\nabla \bl{a}\cdot \bl{v}_Ndx}_{I_{21}}+\un{\int_{\mD}\bl{a}\cdot\nabla \bl{a}\cdot \bl{v}_Ndx}_{I_{22}}.
\]
Noting that $\bl{v}_N\cdot\bl{n}=0$ on $\p\mD$, using integration by parts, one finds
\[
\bali
I_{21}=&\sum_{k,l=1}^3\int_{\mD}(\bl{v}_N)_k\p_{x_k}{a}_l(\bl{v}_N)_jdx=-\sum_{k,l=1}^3\int_{\mD}(\bl{v}_N)_k{a}_l\p_{x_k}(\bl{v}_N)_ldx.
\eali
\]
Using H\"older's inequality, we have
\[
I_{21}\leq C\|\bl{a}\|_{L^\i(\mD)}\|\bl{v}_N\|_{H^1(\mD)}^2.
\]
For the term $I_{22}$, one notices that $\bl{a}$ equals to the Poiseuille flow $\bl{g}^L_{\Phi}$ or $\bl{g}^R_{\Phi}$ in $\mD-\mD_Z$, and thus $\bl{a}\cdot\nabla \bl{a}\equiv 0$ in $\mD-\mD_Z$. This indicates
\[
\begin{split}
|I_{22}|&=\left|\int_{\mD_Z}\bl{a}\cdot\nabla \bl{a}\cdot \bl{v}_Ndx\right|\leq\|\bl{a}\|_{L^3(\mD_Z)}\|\na \bl{a}\|_{L^2(\mD_Z)}\|\bl{v}_N\|_{L^2(\mD)}\leq C_{\mD}\|\bl{a}\|_{H^1(\mD_Z)}^2\|\bl{v}_N\|_{H^1(\mD)}.\\
\end{split}
\]

Finally, it remains to estimate $I_3$. Similarly as $I_{22}$, we also claim that
\[
|I_3|=\left|\int_{\mD_Z}\bl{v}_N\cdot\Dl \bl{a}dx\right|\leq C_Z\|\bl{a}\|_{H^2(\mD_Z)}\|\bl{v}_N\|_{H^1(\mD)}.
\]
Here goes the proof of the claim: By the construction of the Poiseuille flow $\bl{g}_{\Phi}^{L}$, one knows
\[
\int_{\Sigma_L\times(-\i,-Z)}\bl{v}_N\cdot\Dl \bl{a}dx=C\int^{-\i}_{-Z}\int_{\Sigma_L}(\bl{v}_N)_3dx=0.
\]
Here actually we can show that $\int_{\Sigma_L}(\bl{v}_N)_3(x_h,x_3)dx_h$ is independent of $x_3$ by using $\text{div }\bl{v}_N=0$ and $\bl{v}_N\cdot\bl{n}=0$. Then using the compact support of $\bl{v}_N$ we can get the above equality.

Substituting the above estimates for $I_1$--$I_3$, and applying Poincar\'e inequalities in Lemma \ref{P2} and Lemma \ref{TORPOIN}, one derives
\[
{P}(\bl{c}^N)\cdot\bl{c}^N\geq\|\bl{v}_N\|_{H^1(\mD)}\left((C_{\al,\mD}-C\Phi)\|\bl{v}_N\|_{H^1(\mD)}-C_{\mD}(\Phi+\Phi^2)\right),
\]
which guarantees
\[
{P}(\bl{c}^N)\cdot\bl{c}^N\geq0,
\]
provided
\[
\Phi\leq \Phi_0:=C^{-1}C_{\al,\mD}\q\text{and}\q|\bl{c}^N|=\|\bl{v}_N\|_{H^1(\mD)}\geq\f{C_Z(\Phi+\Phi^2)}{C_{\al,\mD}-C\Phi}:=\rho.
\]
Using Lemma \ref{FUNC}, there exists
\be\label{ubound}
\bl{v}_N^*\in\text{span}\left\{\bl{\varphi}_1,\bl{\varphi}_2,...,\bl{\varphi}_N\right\},\q\text{and}\q\|\bl{v}_N^*\|_{H^1(\mD)}\leq\f{C_{\mD}(\Phi+\Phi^2)}{C_{\al,\mD}-C\Phi},
\ee
such that
\begin{equation}\label{app-N}
\begin{split}
&2\int_{\mD}\mathbb{S}\bl{v}_N^*:\mathbb{S}\bl{\phi}_Ndx+\al\int_{\p\mD}(\bl{v}^*_N)_{\mathrm{tan}}\cdot(\bl{\phi}_N)_{\mathrm{tan}}dS+\int_{\mD}\bl{v}_N^*\cdot\nabla \bl{v}_N^*\cdot \bl{\phi}_Ndx+\int_{\mD}\bl{v}_N^*\cdot\nabla \bl{a}\cdot\bl{\phi}_Ndx\\
+&\int_{\mD}\bl{a}\cdot\nabla\bl{v}_N^*\cdot\bl{\phi}_N dx=\int_{\mD}\big(\Dl \bl{a}-\bl{a}\cdot\nabla \bl{a}\big)\cdot\bl{\phi}_Ndx,\q\forall\bl{\phi}_N\in\text{span}\left\{\bl{\varphi}_1,\bl{\varphi}_2,...,\bl{\varphi}_N\right\}.
\end{split}
\end{equation}

The above bound \eqref{ubound} and Rellich-Kondrachov embedding theorem imply the existence of a field $\bl{v}\in \mH_{\sigma}(\mD)$ and a subsequence, which we will always denote by $\bl{v}^*_N$, such that
\begin{equation*}
\bl{v}^*_N\to \bl{v}\quad \text{weakly in $\mH_{\sigma}(\mD)$}
\end{equation*}
and
\begin{equation*}
\bl{v}^*_N\to \bl{v}\quad \text{strongly in $L^2(\mD')$, for all bounded $\mD'\subset\mD$}\,.
\end{equation*}
Therefore, we can pass to the limit in \eqref{app-N} and obtain
\be\label{V}
\begin{split}
&2\int_{\mD}\mathbb{S}\bl{v}:\mathbb{S}\bl{\varphi} dx+\al\int_{\p\mD}\bl{v}_{\mathrm{tan}}\cdot\bl{\varphi}_{\mathrm{tan}}dS+\int_{\mD}\bl{v}\cdot\nabla \bl{v}\cdot\bl{\varphi} dx+\int_{\mD}\bl{v}\cdot\nabla \bl{a}\cdot\bl{\varphi} dx\\
+&\int_{\mD}\bl{a}\cdot\nabla \bl{v}\cdot\bl{\varphi} dx=\int_{\mD}\big(\Dl \bl{a}-\bl{a}\cdot\nabla \bl{a}\big)\cdot\bl{\varphi} dx,\q\text{for any}\q \bl{\varphi}\in \mH_\sigma(\mD).
\end{split}
\ee
Finally, the $H^1$-estimate of $\bl{v}$
\[
\|\bl{v}\|_{H^1(\mD)}\leq C_{\al,\mD}\Phi
\]
follows from \eqref{ubound} and the Fatou lemma for weakly convergent sequences.  This completes the construction of $\bl{v}$.

\subsubsection{Creating the pressure field}\label{SEC332}
\q\ While processing the Galerkin method in the previous subsection, we did nothing with the pressure. This is because all test functions are divergence-free. To find the pressure, we introduce the following lemma, which is a special case of \cite[Theorem 17]{DeRham1960} by De Rham. See also \cite[Proposition 1.1]{Temam1984}.
\begin{lemma}\label{DeRham}
For a given open set $\O\subset\mathbb{R}^3$, let $\bl{\mathcal{F}}$ be a distribution in $\left(C_c^\i(\O;\mathbb{R}^3)\right)'$ which satisfies:
\[
\langle \bl{\mathcal{F}},\bl{\phi}\rangle=0,\q\text{for all}\q \bl{\phi}\in\{\bl{g}\in C_c^\i(\O;\mathbb{R}^3):\,\text{div }\bl{g}=0\}.
\]
Then there exists a distribution $q\in\left(C_c^\i(\O;\mathbb{R})\right)'$ such that
\[
\bl{\mathcal{F}}=\na q.
\]
\end{lemma}

\qed

Let $\bl{v}$ be a weak solution of \eqref{EQU111} constructed in the previous subsection. Using \eqref{V}, one has $\bl{u}=\bl{v}+\bl{a}$ satisfies
\[
\int_{\mD}\nabla \bl{u}\cdot\nabla\bl{\phi}\,dx+\int_{\mD}\bl{u}\cdot\nabla \bl{u}\cdot\bl{\phi}\,dx=0,\q\text{for all}\q \bl{\phi}\in\{\bl{g}\in C_c^\i(\mD;\mathbb{R}^3):\,\text{div }\bl{g}=0\}.
\]
Thus by Lemma \ref{DeRham}, there exists $p\in\left(C_c^\i(\mD;\mathbb{R})\right)'$, such that
\be\label{NS11}
\Dl\bl{u}-\bl{u}\cdot\nabla\bl{u}=\nabla p
\ee
in the sense of distribution.

To derive the regularity of the pressure, one first notices that \eqref{NS11} is equivalent to
\be\label{DEFPi}
\text{div}\big(\nabla\bl{v}-\bl{v}\otimes\bl{v}-\bl{a}\otimes\bl{v}-\bl{v}\otimes\bl{a}\big)+\Dl \bl{a}+\left(\f{\eta(x_3)}{C_{P,R}}+\f{\eta(-x_3)}{C_{P,L}}\right)\Phi\bl{e_3}-\bl{a}\cdot\nabla \bl{a}=\na\Pi,
\ee
with
\be\label{DEFPI}
\Pi=p+\f{\Phi\int_{-\i}^{x_3}\eta(s)ds}{C_{P,R}}-\f{\Phi\int_{-\i}^{-x_3}\eta(s)ds}{C_{P,L}}.
\ee
Here $C_{P,i}$ ($i=L,\,R$) are Poiseuille constants defined in \eqref{POS}, where we had dropped indexes $L$, $R$ there. And $\eta$ is the cut-off function which is given in \eqref{ETA}. By the definition of $\bl{a}$ in \eqref{CONSa}, one has both $\Dl \bl{a}+\left(\f{\eta(x_3)}{C_{P,R}}+\f{\eta(-x_3)}{C_{P,L}}\right)\Phi\bl{e_3}$ and $\bl{a}\cdot\nabla\bl{a}$ are smooth and have compact support. Meanwhile, since $\bl{v}\in H^1(\mD)$ and $\bl{a}$ is uniformly bounded, one deduces
\[
\nabla\bl{v}-\bl{v}\otimes\bl{v}-\bl{a}\otimes\bl{v}-\bl{v}\otimes\bl{a}\in L^2(\mD),
\]
directly by the Sobolev embedding and H\"older's inequality. Therefore one concludes the left hand side of \eqref{DEFPi} belongs to $H^{-1}(\mD)$. Then the following lemma implies $\Pi\in L^2_{\mathrm{loc}}(\overline{\mD})$, which leads to $p\in L^2_{\mathrm{loc}}(\overline{\mD})$ by \eqref{DEFPI}.
\begin{lemma}[See \cite{Temam1984}, Proposition 1.2]\label{LEM312}
Let $\Omega$ be a bounded Lipschitz open set in $\mathbb{R}^{3}$. If a distribution $q$ has all its first derivatives $\p_{x_i} q$, $1\leq i \leq3$, in $H^{-1}(\Omega)$, then $q\in L^{2}(\Omega)$ and
\be\label{EEEE0}
\left\|q-\bar{q}_\O\right\|_{L^{2}(\Omega)}\leq C_\O\|\na q\|_{H^{-1}(\Omega)},
\ee
where $\bar{q}_\O=\f{1}{|\O|}\int_{\O}qdx$. Moreover, if $\Omega$ is any Lipschitz open set in $\mathbb{R}^{3}$, then $q\in L_{\mathrm{loc}}^{2}(\overline{\Omega})$.
\end{lemma}

\qed

This completes the proof of Theorem \ref{THM3.8}.
\section{Uniqueness of the weak solution}\label{sec2}
\q Recall the solution $(\bl{u},{p})$ we constructed in the last section with its flux being $\Phi$. In this part, we will show it is unique if $\Phi$ is sufficiently small.

\subsection{Estimate of the pressure}
\q\ The following lemma states the existence of the solution to problem $\nabla\cdot \bl{V}=f$ in a truncated pipe.
\begin{lemma}\label{LEM2.1}
Let $D=\Sigma\times [0,1]$, $f\in L^2(D)$ with
\[
\int_D fdx=0,
\]
then there exists a vector valued function $\bl{V}:\,D\to\mathbb{R}^3$ belongs to $H^1_0(D)$ such that
\be\label{LEM2.11}
\nabla\cdot \bl{V}=f,\q\text{and}\q\|\nabla \bl{V}\|_{L^2(D)}\leq C\|f\|_{L^2(D)}.
\ee
Here $C>0$ is an absolute constant.
\end{lemma}

See \cite{Bme1, Bme2}, also \cite[Chapter III]{Galdi2011} for detailed proof of this lemma.

\qed

Below, the proposition shows an $L^2$ estimate related to the pressure in the truncated pipe $\O_Z^{+}$ or $\O_Z^{-}$ could be bounded by $L^2$ norm of $\nabla \bl{u}$.
\begin{proposition}\label{P2.4}
Let $(\tilde{\bl{u}},\tilde{p})$ be an alternative weak solution of \eqref{NS} in the pipe $\mD$, subject to the Navier-slip boundary condition \eqref{NBC}. If the total flux
\[
\int_{\mD\cap\{x_3=z\}}\tilde{\bl{u}}(x_h,z)\cdot\boldsymbol{e_3}dx_h=\Phi=\int_{\mD\cap\{x_3=z\}}{\bl{u}}(x_h,z)\cdot\boldsymbol{e_3}dx_h,\q\text{for any } |z|\geq Z,
\]
then the following estimate of $\bl{w}:=\tilde{\bl{u}}-{\bl{u}}$ and the pressure holds
\[
\left|\int_{\O^{\pm}_K}(\tilde{p}-p)w_3dx\right|\leq C_\mD\left(\|\bl{u}\|_{L^4(\O^{\pm}_K)}\|\nabla \bl{w}\|^2_{L^2(\O^{\pm}_K)}+\|\nabla \bl{w}\|_{L^2(\O^{\pm}_K)}^2+\|\nabla \bl{w}\|_{L^2(\O^{\pm}_K)}^3\right),\forall K\geq Z+1,
\]
where $C_\mD>0$ is a constant independent of $K$.
\end{proposition}
\pf During the proof, we cancel the upper index ``$\pm$" of the domain for simplicity. Noticing
\[
\int_{\mD\cap\{x_3=z\}}w_3(x_h,z)dx_h\equiv0,\q\forall |z|\geq Z,
\]
we deduce that
\[
\int_{\O_K}w_3dx=0,\q\forall K\geq Z+1.
\]
Using Lemma \ref{LEM2.1}, one derives the existence of a vector field $V$ satisfying \eqref{LEM2.11} with $f=w_3$. Applying equation \eqref{NS}$_1$, one arrives
\[
\int_{\O_K}(\tilde{p}-{p})w_3dx=-\int_{\O_K}\nabla (\tilde{p}-{p})\cdot \bl{V}dx=\int_{\O_K}\left(\bl{w}\cdot\na \bl{w}+\bl{u}\cdot\nabla \bl{w}+\bl{w}\cdot\nabla \bl{u}-\Dl \bl{w}\right)\cdot \bl{V}dx.
\]
Using integration by parts, one deduces
\[
\int_{\O_K}(\tilde{p}-p)w_3dx=\sum_{i,j=1}^3\int_{\O_K}(\p_iw_j-w_iw_j-u_iw_j-u_jw_i)\p_iV_jdx.
\]
By applying H\"older's inequality and \eqref{LEM2.11} in Lemma \ref{LEM2.1}, one deduces that
\be\label{EP1}
\left|\int_{\O_K}(\tilde{p}-p)w_3dx\right|\leq C\left(\|\na \bl{w}\|_{L^2(\O_K)}+\|\bl{w}\|_{L^4(\O_K)}^2+\|\bl{u}\|_{L^4(\O_K)}\|\bl{w}\|_{L^4(\O_K)}\right)\|w_3\|_{L^2(\O_K)}.
\ee
Since $v_3$ has a zero mean value on each cross section $\Sigma$, and $(\bl{w}-w_3\bl{e_3})$ satisfies
\[
(\bl{w}-w_3\bl{e_3})\cdot\bl{n}=0,\quad\text{for any } x\in\p\mD\cap\p\O_K,
\]
the vector $\bl{w}$ enjoys the Poincar\'e inequality
\be\label{PON3}
\|\bl{w}\|_{L^2(\O_K)}\leq C_\mD  \|\na_h \bl{w}\|_{L^2(\O_K)}.
\ee
Substituting \eqref{PON3} in \eqref{EP1}, also noting the Gagliardo-Nirenberg inequality
\[
\|\bl{w}\|_{L^4(\O_K)}^2\leq C_{\mD}\left(\|\bl{w}\|_{L^2(\O_K)}^{1/2}\|\na \bl{w}\|_{L^2(\O_K)}^{3/2}+\|\bl{w}\|_{L^2(\O_K)}^{2}\right),
\]
one concludes
\[
\left|\int_{\O_K}(\tilde{p}-p)w_3dx\right|\leq C_\mD\left(\|\bl{u}\|_{L^4(\O_K)}\|\nabla \bl{w}\|^2_{L^2(\O_K)}+\|\nabla \bl{w}\|_{L^2(\O_K)}^2+\|\nabla \bl{w}\|_{L^2(\O_K)}^3\right).
\]

\qed

\subsection{Main estimates}
\q\ Subtracting the equation of $\bl{u}$ from the equation of $\tilde{\bl{u}}$, one finds
\be\label{SUBT}
\bl{w}\cdot\nabla \bl{w}+\bl{u}\cdot\nabla \bl{w}+\bl{w}\cdot\nabla \bl{u}+\nabla(\tilde{p}-p)-\Dl \bl{w}=0.
\ee
Multiplying $\bl{w}$ on both sides of \eqref{SUBT}, and integrating on $\mathcal{D}_\zeta$, one derives
\be\label{Maint0}
\int_{\mathcal{D}_\zeta}\bl{w}\cdot\Dl \bl{w}dx=\int_{\mathcal{D}_\zeta}\bl{w}\big(\bl{w}\cdot\nabla \bl{w}+\bl{u}\cdot\nabla \bl{w}+\bl{w}\cdot\nabla \bl{u}+\nabla (\tilde{p}-p)\big)dx.
\ee
Using the divergence-free property and the Navier-slip boundary condition of $\bl{u}$ and $\tilde{\bl{u}}$, one deduces
\be\label{Maint1}
\begin{split}
\int_{\mathcal{D}_\zeta}\bl{w}\cdot\Dl \bl{w}dx&=\int_{\mD_\zeta}w_i\p_{x_j}(\p_{x_j}w_i+\p_{x_i}w_j)dx\\
&=-\sum_{i,j=1}^3\int_{\mD_\zeta}\p_{x_j}w_i(\p_{x_j}w_i+\p_{x_i}w_j)dx+\sum_{i,j=1}^3\int_{\p\mD_\zeta}w_in_j(\p_{x_j}w_i+\p_{x_i}w_j)dx\\
&=-2\int_{\mD_\zeta}|\mathbb{S} \bl{w}|^2dx-\al\int_{\p\mD_\zeta\cap\p\mD}\left(|w_{\tau_1}|^2+|w_{\tau_2}|^2\right)dS\\
&\hskip 1cm+\sum_{i=1}^3\int_{\Sigma_R\times\{x_3=\zeta\}}w_i(\p_{x_3}w_i+\p_{x_i}w_3)dx_h-\sum_{i=1}^3\int_{\Sigma_L\times\{x_3=-\zeta\}}w_i(\p_{x_3}w_i+\p_{x_i}w_3)dx_h.
\end{split}
\ee
Here $\bl{n}=(n_1,n_2,n_3)$ is the unit outer normal vector on $\p\mD$.

On the other hand, using integration by parts, we may derive alternatively:
\be\label{Maint20}
\begin{split}
\int_{\mathcal{D}_\zeta}\bl{w}\cdot\Dl \bl{w}dx=-\int_{\mD_\zeta}|\nabla \bl{w}|^2dx+\un{\f{1}{2}\int_{\p\mD_\zeta}\nabla |\bl{w}|^2\cdot\boldsymbol{n}dS}_{T_1},
\end{split}
\ee
where
\[
\begin{split}
T_1=&\un{\f{1}{2}\int_{\p\mD_\zeta\cap\p\mD_M}\nabla |\bl{w}|^2\cdot\boldsymbol{n}dS}_{T_{11}}+\un{\f{1}{2}\int_{\p\mD_\zeta\cap\p\mD_L}\nabla |\bl{w}|^2\cdot\boldsymbol{n}dS}_{T_{12}}+\un{\f{1}{2}\int_{\p\mD_\zeta\cap\p\mD_R}\nabla |\bl{w}|^2\cdot\boldsymbol{n}dS}_{T_{13}}\\
&+\f{1}{2}\left(\int_{\mD\cap\{x_3=\zeta\}}\p_{x_3}|\bl{w}|^2(x_h,\zeta)dx_h-\int_{\mD\cap\{x_3=-\zeta\}}\p_{x_3}|\bl{w}|^2(x_h,-\zeta)dx_h\right).
\end{split}
\]
To bound the term $T_{11}$, we will apply the local orthogonal curvilinear coordinates on $\p\mD$. Thus we split $\p\mD_{\zeta}\cap\p\mD_M$ into finitely many pieces:
\[
\p\mD_{\zeta}\cap\p{\mD_M}=\bigcup_{i=1}^NV_i,
\]
 and in each piece $V_i$, there exists an orthogonal curvilinear frame $\{\bl{\tau_1}^i,\bl{\tau_2}^i,\bl{n}^i\}$ such that
\[
\nabla |\bl{w}|^2=\p_{\bl{\tau_1}^i}|\bl{w}|^2\boldsymbol{\tau_1}^i+\p_{\bl{\tau_2}^i}|\bl{w}|^2\boldsymbol{\bl{\tau_2}}^i+\p_{\bl{n}^i}|\bl{w}|^2\boldsymbol{n}^i,\q\text{on }V_i.
\]
Using \eqref{NBCM}, one derives
\be\label{T11}
\begin{split}
|T_{11}|\leq&\sum_{i=1}^N\int_{V_i}\left|w_{\tau_1^i}\left(\al-\kappa_1^i(x)\right)w_{\tau_1^i}\right|dS+\sum_{i=1}^N\int_{V_i}\left|w_{\tau_2^i}\left(\al-\kappa_2^i(x)\right)w_{\tau_2^i}\right|dS\\
\leq &C_{\al,\mD}\int_{\p\mD_{\zeta}\cap\p{\mD_M}}|\bl{w}_{\mathrm{tan}}|^2dS.
\end{split}
\ee
Here $C_{\al,\mD}>0$ is a constant depending only on the friction ratio $\al$ and the domain $\mD$. The existence of this constant $C_{\al,\mD}$ follows from the boundedness of the principal curvature of $\p\mD$. See Proposition \ref{Prop22} for details.

Noting that $\mD_L$ is a straight pipe, one can find the global natural coordinates $\{\bl{\tau_1},\bl{e_3},\bl{n}\}$ of $\p\mD\cap\p{\mD_L}$, where $\bl{\tau_1}$ and $\bl{n}$ are unit tangent and unit outer normal vector of $\p\Sigma_L$ in the $x_1Ox_2$ plane, while $\bl{e_3}$ is the Euclidean coordinate vector in $x_3$-direction. In this case, one writes
\[
\nabla |\bl{w}|^2=\p_{\bl{\tau_1}}|\bl{w}|^2\boldsymbol{\tau_1}+\p_{x_3}|\bl{w}|^2\boldsymbol{e_3}+\p_{\bl{n}}|\bl{w}|^2\boldsymbol{n},\q\text{on }\p\mD_{\zeta}\cap\p{\mD_L}.
\]
Using \eqref{NBCM1}, one has $T_{12}$ satisfies
\be\label{T12}
|T_{12}|\leq\left|\int_{\p\mD_{\zeta}\cap\p{\mD_L}}\big(w_{\tau}\left(\al-\kappa_1(x)\right)w_{\tau}+\al|w_3|^2\big)dS\right|\leq C_{\al,\mD}\int_{\p\mD_{\zeta}\cap\p{\mD_L}}|\bl{w}_{\mathrm{tan}}|^2dS.\\
\ee
Similarly as \eqref{T12}, one derives
\be\label{T13}
|T_{13}|\leq C_{\al,\mD}\int_{\p\mD_{\zeta}\cap\p{\mD_R}}|\bl{w}_{\mathrm{tan}}|^2dS.\\
\ee
Substituting \eqref{T11}--\eqref{T13} in \eqref{Maint20}, one concludes that
\be\label{Maint2}
\int_{\mathcal{D}_\zeta}\bl{w}\cdot\Dl \bl{w}dx\leq-\int_{\mD_\zeta}|\nabla \bl{w}|^2dx+C_{\al,\mD}\int_{\p\mD_\zeta\cap\p\mD}|\bl{w}_{\mathrm{tan}}|^2dS+C\int_{\mD\cap\{x_3=\pm\zeta\}}|\bl{w}||\nabla \bl{w}|dx_h.
\ee
Now we focus on the right hand side of \eqref{Maint0}. Applying integration by parts, one derives
\be\label{Maint3}
\begin{split}
\int_{\mathcal{D}_\zeta}\bl{w}\big(\bl{w}\cdot\nabla \bl{w}+\nabla (\tilde{p}-p)\big)dx=&\int_{\mD\cap\{x_3=\zeta\}}w_3\left(\f{1}{2}|\bl{w}|^2+(\tilde{p}-p)\right)dx\\
&-\int_{\mD\cap\{x_3=-\zeta\}}w_3\left(\f{1}{2}|\bl{w}|^2+(\tilde{p}-p)\right)dx.
\end{split}
\ee
Applying H\"older's inequality, noting that $\bl{u}=\bl{v}+\bl{a}$, where $\bl{a}$ is the profile vector defined in Section \ref{SEC3.2}, while $\bl{v}$ is the $H^1$-weak solution constructed in Section \ref{SEC3627}, one has
\be\label{Maint5}
\begin{split}
\left|\int_{\mathcal{D}_\zeta}\big(\bl{w}\cdot\nabla \bl{u}\cdot \bl{w}+\bl{u}\cdot\na\bl{w}\cdot\bl{w}\big)dx\right|\leq&\, \|\nabla\bl{v}\|_{L^2(\mD_\zeta)}\|\bl{w}\|_{L^4(\mD_\zeta)}^2+\|\bl{v}\|_{L^4(\mD_\zeta)}\|\na\bl{w}\|_{L^2(\mD_\zeta)}\|\bl{w}\|_{L^4(\mD_\zeta)}\\
&+\|\na\bl{a}\|_{L^\i(\mD_\zeta)}\|\bl{w}\|_{L^2(\mD_\zeta)}^2+\|\bl{a}\|_{L^\i(\mD_\zeta)}\|\na\bl{w}\|_{L^2(\mD_\zeta)}\|\bl{w}\|_{L^2(\mD_\zeta)}\\
\leq&\, C_{\mD}\left(\|\bl{v}\|_{H^{1}(\mD_\zeta)}+\|\bl{a}\|_{W^{1,\i}(\mD_\zeta)}\right)\int_{\mD_\zeta}|\nabla \bl{w}|^2dx\\
\leq&\,C_{\al,\mD}\Phi\int_{\mD_\zeta}|\nabla \bl{w}|^2dx.
\end{split}
\ee
Here in the second inequality, we have applied the Gagliardo-Nirenberg inequality and the Poincar\'e inequality \eqref{TORPIPEPOIN} in Lemma \ref{TORPOIN}, which indicate
\[
\|\bl{w}\|_{L^4(\mD_\zeta)}\leq C_{\mD}\left(\|\bl{w}\|^{1/4}_{L^2(\mD_\zeta)}\|\na\bl{w}\|^{3/4}_{L^2(\mD_\zeta)}+\|\bl{w}\|_{L^2(\mD_\zeta)}\right)\leq C_\mD\left(\int_{\mD_\zeta}|\nabla \bl{w}|^2dx\right)^{1/2}.
\]
 Meanwhile, the third inequality in \eqref{Maint5} is guaranteed by \eqref{EOFV627} and \eqref{EOFAA6277}.

Therefore, by calculating
\[
\eqref{Maint1}\times C_{\al,\mD}+\eqref{Maint2}\times\al,
\]
we derive
\be\label{Maint4}
\int_{\mD_\zeta}\bl{w}\cdot\Dl \bl{w}dx\leq-2C_{\al,\mD}\int_{\mD_\zeta}|\mathbb{S}\bl{w}|^2dx-\al\int_{\mD_\zeta}|\nabla \bl{w}|^2dx+C\int_{\mD\cap\{x_3=\pm\zeta\}}|\bl{w}||\na \bl{w}|dx_h.
\ee
Substituting \eqref{Maint3}, \eqref{Maint5} and \eqref{Maint4} in \eqref{Maint0}, one arrives
\[
\begin{split}
\left(1-C_{\al,\mD}\Phi\right)\int_{\mD_\zeta}|\nabla \bl{w}|^2dx\leq& \f{C_{\al,\mD}}{\al}\left(\int_{\mD\cap\{x_3=\pm\zeta\}}|\bl{w}|(|\na \bl{w}|+|\bl{w}|^2)dx_h\right.\\
&\left.-\int_{\mD\cap\{x_3=\zeta\}}w_3\left(\tilde{p}-p\right)dx_h+\int_{\mD\cap\{x_3=-\zeta\}}w_3\left(\tilde{p}-p\right)dx_h\right).
\end{split}
\]
Now one concludes that if $\Phi<<1$ being small enough such that $C_{\al,\mD}\Phi<\f{1}{2}$
\[
\begin{split}
\int_{\mD_\zeta}|\nabla \bl{w}|^2dx\leq \f{C_{\al}}{\al\left(1-C_{\al,\mD}\Phi\right)}&\left(\int_{\mD\cap\{x_3=\pm\zeta\}}|\bl{w}|(|\na \bl{w}|+|\bl{w}|^2)dx_h\right.\\
&\left.-\int_{\mD\cap\{x_3=\zeta\}}w_3\left(\tilde{p}-p\right)dx_h+\int_{\mD\cap\{x_3=-\zeta\}}w_3\left(\tilde{p}-p\right)dx_h\right).
\end{split}
\]
Therefore,  one derives the following estimate by integrating with $\zeta$ on $[K-1,K]$, where $K\geq Z+1$:
\be\label{ET+0}
\int_{K-1}^K\int_{\mathcal{D}_\zeta}|\nabla \bl{w}|^2dxd\zeta\leq C_{\al,\mD}\Bigg(\int_{\O_K^+\cup\O_K^-}|\bl{w}|(|\na \bl{w}|+|\bl{w}|^2)dx+\Big|\int_{\O_K^+\cup\O_K^-}w_3\left(\tilde{p}-p\right)dx\Big|\Bigg).
\ee
Here $C_{\al,\mD}>0$ is a constant. Now we only handle integrations on $\O_K^+$ since the cases of $\O_K^-$ are similar. Using the Cauchy-Schwarz inequality and the Poincar\'e inequality Lemma \ref{P2}, one has
\be\label{ET+1}
\int_{\O_K^+}|\bl{w}||\na \bl{w}|dx\leq\|\bl{w}\|_{L^2(\O_K^+)}\|\nabla \bl{w}\|_{L^2(\O_K^+)}\leq C\|\nabla \bl{w}\|^2_{L^2(\O_K^+)}.
\ee
Moreover, by H\"older's inequality and the Gagliardo-Nirenberg inequality, one writes
\[
\int_{\O_K^+}|\bl{w}|^3dx\leq C_\mD\left(\|\bl{w}\|_{L^2(\O_K^+)}^{3/2}\|\na \bl{w}\|_{L^2(\O_K^+)}^{3/2}+\|\bl{w}\|_{L^2(\O_K^+)}^{3}\right),
\]
which follows by the Poincar\'e inequality that
\[
\int_{\O_K^+}|\bl{w}|^3dx\leq C\|\na \bl{w}\|_{L^2(\O_K^+)}^{3}.
\]
Recalling Proposition \ref{P2.4}, one arrives at
\be\label{ET+2}
\left|\int_{\O_K^+}w_3\left(\tilde{p}-p\right)dx\right|\leq C\left(\|\bl{u}\|_{L^4(\O_K^+)}\|\nabla \bl{w}\|^2_{L^2(\O_K^+)}+\|\nabla \bl{w}\|_{L^2(\O_K^+)}^2+\|\nabla \bl{w}\|_{L^2(\O_K^+)}^3\right).
\ee
Substituting \eqref{ET+1}--\eqref{ET+2}, together with their related inequality on domain $\O_K^{-}$, in \eqref{ET+0}, one concludes
\be\label{FEST}
\int_{K-1}^K\int_{\mathcal{D}_\zeta}|\nabla \bl{w}|^2dxd\zeta\leq C_{\al,\mD}\left(\|\nabla \bl{w}\|_{L^2(\O_K^+\cup\O_K^-)}^2+\|\nabla \bl{w}\|_{L^2(\O_K^+\cup\O_K^-)}^3\right).
\ee
\subsection{End of proof}

\q\ Finally, by defining
\[
Y(K):=\int_{K-1}^K\int_{\mathcal{D}_\zeta}|\nabla \bl{w}|^2dxd\zeta,
\]
\eqref{FEST} indicates
\[
Y(K)\leq C_{\al,\mD}\left(Y'(K)+\left(Y'(K)\right)^{3/2}\right),\q\forall K\geq Z+1.
\]
By Lemma \ref{LEM2.3}, we derive
\[
\liminf_{\zeta\to\infty}K^{-3}Y(K)>0,
\]
that is, there exists $C_0>0$ such that
\[
\int_{K-1}^K\int_{\mathcal{D}_\zeta}|\nabla \bl{w}|^2dxd\zeta\geq C_0K^3.
\]
However, this leads to a paradox with the condition \eqref{GROWC}. Thus $Y(K)\equiv0$ for all $K\geq Z+1$, which proves $\bl{u}\equiv \tilde{\bl{u}}$.

\qed



\section{Regularity and decay estimates of the weak solution}\label{SEC5}
\q\ In this section, we will show the weak solution, which is proved to be unique in Section \ref{sec2}, is smooth, and it decays exponentially to Poiseuille flows $\bl{g}_\Phi^R$ and $\bl{g}_\Phi^L$ as $x_3\to\pm\i$, respectively. Recall
\be\label{DEFVPI}
\bl{v}=\bl{u}-\bl{a},\q \Pi=p+\f{\Phi\int_{-\i}^{x_3}\eta(s)ds}{C_{P,R}}-\f{\Phi\int_{-\i}^{-x_3}\eta(s)ds}{C_{P,L}}.
\ee
The route of the proof follows that:
\begin{itemize}
\setlength{\itemsep}{-2 pt}
\item[(i).] Global $W^{1,3}$ and $H^2$ estimate of $\bl{v}$, together with the global $L^2$-estimate of $\na\Pi$;
\item[(ii).] Higher-order regularity of $(\bl{v},\Pi)$;
\item[(iii).] $H^1$-exponential decay estimate of $\bl{v}$;
\item[(iv).] The exponential decay for higher-order norms of $(\bl{v},\Pi)$.
\end{itemize}

\subsection{Global $\bl{H^2}$ estimate of the solution}
\q\ In this subsection, we will show the weak solution constructed in Section \ref{SECE} is strong. Our strategy is treating the Navier-Stokes system \eqref{EQU} as the following linear Stokes equations:
\be\label{EQU1}
\left\{
\begin{array}{*{2}{ll}}
-\Dl \bl{v}+\na \Pi=\mathrm{div}\bl{F}+\bl{f}, \q\na\cdot \bl{v}=0,\q\q&\text{in }\mD;\\[2mm]
2(\mathbb{S}\bl{v}\cdot\bl{n})_{\mathrm{tan}}+\al \bl{v}_{\mathrm{tan}}=0,\q \bl{v}\cdot\bl{n}=0, \q\q&\text{on }\p\mD,
\end{array}
\right.
\ee
where
\[
\begin{split}
\bl{F}:&=-\bl{v}\otimes\bl{v}-\bl{a}\otimes\bl{v}-\bl{v}\otimes\bl{a};\\
\bl{f}:&=\Dl \bl{a}+\Phi\left(\f{\eta(x_3)}{C_{P,R}}+\f{\eta(-x_3)}{C_{P,L}}\right)\bl{e_3}-\bl{a}\cdot\nabla \bl{a},
\end{split}
\]
and then applying the bootstrapping method. Noting that $\bl{v}\in H^{1}(\mD)$, the Sobolev imbedding indicates $\bl{v}\in L^2(\mD)\cap L^6(\mD)$, which indicates
\[
\bl{v}\otimes\bl{v}\in L^3(\mD).
\]
Moreover, noting that $\bl{a}$ is smooth and uniformly bounded in $\mD$, and both $\Dl \bl{a}+\Phi\left(\f{\eta(x_3)}{C_{P,R}}+\f{\eta(-x_3)}{C_{P,L}}\right)\bl{e_3}$ and $\bl{a}\cdot\nabla\bl{a}$ have compact support, one concludes
\be\label{L3EST}
\bl{F}\in L^3(\mD),\q\bl{f}\in C_c^\i(\overline{\mD}).
\ee
Here goes the main result of this subsection:
\begin{proposition}\label{PROP51}
Let $(\bl{u},p)$ be the weak solution to \eqref{GL1}--\eqref{GL3}, and $(\bl{v}, \Pi)$ is defined as \eqref{DEFVPI}. Then
\begin{equation}\label{reg0}
(\bl{v}\,,\na\Pi)\in H^2(\mD)\times L^2(\mD),
\end{equation}
which satisfies
\be\label{hesti}
\|\bl{v}\|_{H^{2}(\mD)}+\|\nabla \Pi\|_{L^2(\mD)}\leq C_{\al,\mD}\Phi.
\ee
\end{proposition}

\pf The proof consists of two parts. First we show $\bl{v}\in W^{1,3}(\mD)$ by applying the regularity results in \eqref{L3EST}. This will conclude that
\be\label{H1EST}
\mathrm{div}\,\bl{F}\in L^{2}(\mD).
\ee
Then based on \eqref{H1EST}, we can obtain \eqref{reg0}.

Now we split the problem \eqref{EQU1} into a sequence of problems on bounded domains so that Lemma \ref{LEM52} and Lemma \ref{LEM59} are valid for each one of them. To do this, we denote
\[
\mD=\bigcup_{k\in\mathbb{Z}}\mfD_k,\q\text{where}\q\mfD_k:=\mD\cap\left\{x\in\mathbb{R}^3:\,\left(\f{3k}{2}-1\right)Z\leq x_3\leq\left(\f{3k}{2}+1\right)Z\right\},
\]
and the related cut-off function $\psi_k=\psi\left(x_3-\f{3kZ}{2}\right)$ which satisfies
\[
\left\{
\begin{array}{*{2}{ll}}
\mathrm{supp}\,\psi\subset\left[-9Z/10\,,9Z/10\right];\\
\psi\equiv 1,&\q\text{in }\left[-{4Z}/{5}\,,{4Z}/{5}\right];\\
0\leq\psi\leq1,&\q\text{in }\left[-Z\,,Z\right];\\
|\psi^{(m)}|\leq C/Z^m\leq C,&\q\text{for }m=1,2. \\
\end{array}
\right.
\]
\begin{remark}
According to the splitting and constructions above, one notices that the ``bubble part" in $\mD$ is totally contained in $\mfD_0$, and $\psi_k'$, for each $k\in\mathbb{Z}$, are totally supported away from the ``bubble part" of $\mD$. Moreover, any point in $\mD$ is contained in at most two neighboring $\mfD_k$, while the union of sets
\[
\mfD_k':=\{x_3\in\mfD_k :\, \psi_k(x_3)=1\},\q k\in\mathbb{Z}
\]
covers $\mD$.
\end{remark}

\qed

Multiplying the linearized equations \eqref{EQU1}$_1$ with $\psi_k$, then we convert the problem \eqref{EQU1} to a related problem in domain $\mfD_k$, with $k\in\mathbb{Z}$.
\be\label{cutEq}
\left\{
\begin{array}{*{2}{ll}}
\begin{split}
-\Dl(\psi_k\bl{v})+\nabla\left(\psi_k(\Pi-\overline{\Pi}_{\mfD_k})\right)=&\mathrm{div}\,\left(\psi_k\bl{F}\right)+\psi_k\bl{f}+\left(\Pi-\overline{\Pi}_{\mfD_k}\right)\psi_k'\bl{e_3}\\
&-2\psi_k'\p_{x_3}\bl{v}-\psi_k''v_3\bl{e_3}-\psi_k'\bl{F}\cdot\bl{e_3}
\end{split},&\q\q\text{in }\mfD_k;\\[2mm]
\nabla\cdot\left(\psi_k\bl{v}\right)=v_3\psi_k',&\q\q\text{in }\mfD_k;\\[2mm]
2(\mathbb{S}(\psi_k\bl{v})\cdot\bl{n})_{\mathrm{tan}}+\al (\psi_k\bl{v})_{\mathrm{tan}}=0,\q (\psi_k\bl{v})\cdot\bl{n}=0, &\q\q\text{on }\p\mfD_k\cap\p\mD.
\end{array}
\right.
\ee
Here the first two lines in \eqref{cutEq} follow from direct calculations, thus we only give some explanation of the boundary condition in \eqref{cutEq}$_3$. According to the construction of the cut-off function $\psi_k$, one knows that for any $k\in\mathbb{Z}$, $\psi_k'$ is supported away from the ``bubble part" of $\mD$. This indicates the Navier-slip boundary condition \eqref{EQU1}$_2$ enjoys the following form  in the orthogonal curvilinear coordinates on the boundary:
\be\label{BCCO}
\left\{
\begin{aligned}
&\p_{\boldsymbol{n}}v_{\tau_1}=\left(\kappa_1(x)-\al\right)v_{\tau_1},\\
&\p_{\boldsymbol{n}}v_3=-\al v_3,\\
&v_n=0.\\
\end{aligned}
\right.\q\q\q\text{on}\q\p\mD\cap(\p\mD_R\cup\p\mD_L).
\ee
See Remark \ref{RMK2.2} for details. Noting that the normal vector $\bl{n}$ depends only on $x_h$ in the ``straight part" of $\mD$, while $\psi_k$ depending only on $x_3$, one deduces the boundary condition of $\psi_k\bl{v}$ shares the form as \eqref{BCCO}. This indicates the validity of \eqref{cutEq}$_3$.

However, Lemma \ref{LEM52} is not legal for $\psi_k\bl{v}$ at the moment, because $\psi_k\bl{v}$ is not divergence-free, also it does not lie on a smooth domain. Nevertheless, let $\tilde{\mfD}_k$ be a bounded smooth domain which contains $\mfD_k$, with its boundary $\p\tilde{\mfD}_k\supset\p\mfD_k\cap\p\mD$. Guaranteed by the definition of $\mfD_k$, one chooses each $\tilde{\mfD}_k$ with $k>0$ to be congruent to $\tilde{\mfD}_1$, and every $\tilde{\mfD}_k$ with $k<0$ to be congruent to $\tilde{\mfD}_{-1}$.
\begin{figure}[H]
\centering
\includegraphics[scale=0.5]{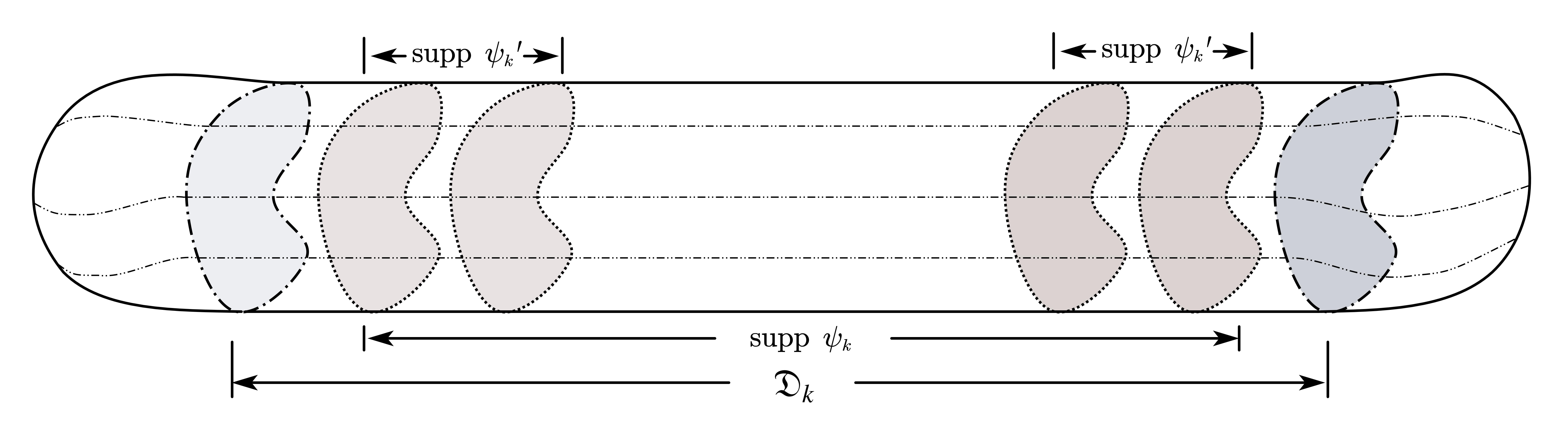}
\caption{Truncated smooth capsule $\tilde{\mfD}_k$.}
\end{figure}

In order to eliminate the divergence part of $\psi_k\bl{v}$, we introduce auxiliary functions $\xi_k$ which satisfy
\be\label{AUXHAR}
\left\{
\begin{aligned}
&\Dl\xi_k=v_3\psi_k',\q\text{in }\tilde{\mfD}_k,\\
&\f{\p\xi_k}{\p\bl{n}}=0,\hskip .9cm\text{on }\p\tilde{\mfD}_k,\\
&\int_{\tilde{\mfD}_k}\xi_k(x)dx=0.
\end{aligned}
\right.
\ee

In the below we denote $\bl{\mathfrak{u}}_k:=\psi_k\bl{v}-\nabla\xi_k$ for convenience. From \eqref{cutEq} and \eqref{AUXHAR}, we know $\bl{\mathfrak{u}}_k$ solves
\be\label{cutEq1}
\left\{
\begin{array}{*{2}{ll}}
-\Dl\bl{\mathfrak{u}}_k+\nabla\left(\psi_k(\Pi-\overline{\Pi}_{\mfD_k})\right)=\mathrm{div}\,\bl{F}_k+\bl{f}_k,&\q\q\text{in }\tilde{\mfD}_k;\\[2mm]
\nabla\cdot\bl{\mathfrak{u}}_k=0,&\q\q\text{in }\tilde{\mfD}_k;\\[2mm]
2(\mathbb{S}\bl{\mathfrak{u}}_k\cdot\bl{n})_{\mathrm{tan}}+\al (\bl{\mathfrak{u}}_k)_{\mathrm{tan}}=\bl{h}_k,\q \bl{\mathfrak{u}}_k\cdot\bl{n}=0, &\q\q\text{on }\p\tilde{\mfD}_k.
\end{array}
\right.
\ee
Here
\[
\begin{split}
\bl{F}_k&:=\psi_k\bl{F}+\nabla^2\xi_k;\\[2mm]
\bl{f}_k&:=\psi_k\bl{f}+\left(\Pi-\overline{\Pi}_{\mfD_k}\right)\psi_k'\bl{e_3}-2\psi_k'\p_{x_3}\bl{v}-\psi_k''v_3\bl{e_3}-\psi_k'\bl{F}\cdot\bl{e_3};\\[2mm]
\bl{h}_k&:=-2\left(\left(\mathbb{S}\nabla\xi_k\right)\cdot\bl{n}\right)_{\mathrm{tan}}-\al\left(\na\xi_k\right)_{\mathrm{tan}}.
\end{split}
\]
Now we are ready to show the regularity estimate of the quantities $\bl{F}_k$, $\bl{f}_k$ and $\bl{h}_k$ above:
\begin{lemma}\label{LEMCOR}
The following estimate of $\bl{F}_k$, $\bl{f}_k$ and $\bl{h}_k$ holds
\be\label{EEEST}
\begin{split}
\|\bl{F}_k\|_{L^3(\tilde{\mfD}_k)}+\|\bl{f}_k\|_{L^{\frac32}(\tilde{\mfD}_k)}+\|\bl{h}_k\|_{W^{-\f{1}{3},3}(\p\tilde{\mfD}_k)}\leq C_{\mD}\|\bl{v}\|_{H^1(\mfD_k)}\left(1+\|\bl{v}\|_{H^1(\mfD_k)}\right)+C\Phi\chi_{|k|\leq1}.\\
\end{split}
\ee
Here the constant $C_{\mD}$ is uniform with $k$, and $\chi_{|k|\leq1}$ is the characteristic function defined by
\[
\chi_{|k|\leq 1}=\left\{
\begin{array}{*{2}{ll}}
1,&\q\q\text{if }\q k\in\{0,1,-1\};\\[2mm]
0,&\q\q\text{if }\q k\in\mathbb{Z}-\{0,1,-1\}.
\end{array}
\right.
\]
\end{lemma}
\pf Noting that the support of $\psi_k$ is uniformly bounded, the estimates of $\bl{F}_k$ and $\bl{h}_k$ in \eqref{EEEST} follows directly from the classical elliptic estimate of system \eqref{AUXHAR}, which states
\be\label{EEEE-1}
\|\xi_k\|_{W^{2,3}(\tilde{\mfD}_k)}\leq C_{\mD}\|\bl{v}\|_{L^3(\mfD_k)}\leq C_{\mD}\|\bl{v}\|_{H^1(\mfD_k)},
\ee
and the trace theorem of Sobolev functions, which follows that
\be\label{EEEE00}
\|\bl{h}_k\|_{W^{-\f{1}{3},3}(\p\tilde{\mfD}_k)}\leq C_{\al,\mD}\|\xi_k\|_{W^{2,3}(\tilde{\mfD}_k)}.
\ee
For the term $\bl{f}_k$, we only derive the estimate of the pressure term since the rest are transparent. Using Lemma \ref{LEM312} and H\"older's inequality, one deduces
\be\label{EEEE1}
\left\|\left(\Pi-\overline{\Pi}_{\mfD_k}\right)\psi_k'\right\|_{L^\frac32(\tilde{\mfD}_k)}\leq C\|\Pi-\overline{\Pi}_{\mfD_k}\|_{L^2(\mfD_k)}\leq C_1\|\na\Pi\|_{H^{-1}(\mfD_k)}.
\ee
Notice that, by the definition of $\mfD_k$ and $\tilde{\mfD}_k$, one has each $\mfD_k$ ($k\in\mathbb{Z}$) is congruent to an element in $\{\mfD_{-1},\,\mfD_0,\,\mfD_1\}$, while every $\tilde{\mfD}_k$  ($k\in\mathbb{Z}$) is congruent to an element in $\{\tilde{\mfD}_{-1},\,\tilde{\mfD}_0,\,\tilde{\mfD}_1\}$. Thus constants in estimates \eqref{EEEE-1}, \eqref{EEEE00} and \eqref{EEEE1} above could be chosen uniformly with respect to $k\in\mathbb{Z}$. Finally, by equation

\[
\nabla\Pi=\text{div}\big(\nabla\bl{v}-\bl{v}\otimes\bl{v}-\bl{a}\otimes\bl{v}-\bl{v}\otimes\bl{a}\big)+\Dl \bl{a}+\left(\f{\eta(x_3)}{C_{P,R}}+\f{\eta(-x_3)}{C_{P,L}}\right)\Phi\bl{e_3}-\bl{a}\cdot\nabla \bl{a},
\]
with both $\Dl \bl{a}+\left(\f{\eta(x_3)}{C_{P,R}}+\f{\eta(-x_3)}{C_{P,L}}\right)\Phi\bl{e_3}$ and $\bl{a}\cdot\nabla \bl{a}$ vanish in $\mfD_k$ with $|k|\geq 2$, one concludes from \eqref{EEEE1} that
\be\label{EEEE11}
\begin{split}
\left\|\left(\Pi-\overline{\Pi}_{\mfD_k}\right)\psi_k'\right\|_{L^\frac32(\tilde{\mfD}_k)}&\leq C\left(\|\na\bl{v}\|_{L^2(\mfD_k)}+\|\bl{v}\|_{L^4(\mfD_k)}^2+\Phi\|\bl{v}\|_{L^2(\mfD_k)}\right)+C\Phi\chi_{|k|\leq1}\\
&\leq C_{\al,\mD}\|\bl{v}\|_{H^1(\mfD_k)}\left(1+\|\bl{v}\|_{H^1(\mfD_k)}\right)+C\Phi\chi_{|k|\leq1}.
\end{split}
\ee
Here we have applied the Sobolev imbedding theorem and interpolations of $L^p$ spaces. This completes the proof of Lemma \ref{LEMCOR}.

\qed

Therefore one concludes the following intermediate $W^{1,3}$ estimate of $\bl{v}$
\be\label{ESTW13}
\|\boldsymbol{v}\|_{{W}^{1, 3}(\mD)}\leq C_{\al,\mD}\Phi
\ee
by combining \eqref{ESTLEM27} in Lemma \ref{LEM52} and \eqref{EEEST} in Lemma \ref{LEMCOR}, then summing up with $k\in\mathbb{Z}$. The details are as follows
\[
\begin{split}
\|\bl{v}\|^3_{W^{1,3}(\mD)}&\leq C\sum_{k\in\mathbb{Z}}\|\psi_k\bl{v}\|^3_{W^{1,3}(\tilde{\mfD}_k)}\\
&\leq C\sum_{k\in\mathbb{Z}}\left(\|\bl{\mathfrak{u}}_k\|^3_{W^{1,3}(\tilde{\mfD}_k)}+\|\xi_k\|^3_{W^{2,3}(\tilde{\mfD}_k)}\right)\\
&\leq C_{\al,\mD}\sum_{k\in\mathbb{Z}}\left(\|\bl{v}\|^3_{H^1(\mfD_k)}\left(1+\|\bl{v}\|^3_{H^1(\mfD_k)}\right)+\Phi^3\chi_{|k|\leq1}\right).
\end{split}
\]
Here we have applied the fact that any point in $\mD$ is contained in at most two neighboring ${\mfD}_k$. Noting that $\Phi\leq\Phi_0=\Phi_0(\al,\mD)$, this further indicates that
\[
\|\bl{v}\|^3_{W^{1,3}(\mD)}\leq C_{\al,\mD}\left((\Phi+\Phi^4)\sum_{k\in\mathbb{Z}}\|\bl{v}\|^2_{H^1(\mfD_k)}+\Phi^3\right)\leq C_{\al,\mD}\Phi^3.
\]

Moreover, \eqref{ESTW13} and H\"older's inequality indicate that
\[
\|\bl{v}\cdot\na\bl{v}\|_{L^2(\mD)}\leq \|\bl{v}\|_{L^6(\mD)}\|\na\bl{v}\|_{L^3(\mD)}<\i,
\]
which further implies
\[
\mathrm{div}\,\bl{F}+\bl{f}\in L^2(\mD).
\]
Similarly as \eqref{EEEST}, now we can deduce
\be\label{EEEST2}
\begin{split}
\left\|\mathrm{div}\,\bl{F}_k+\bl{f}_k\right\|_{L^2(\tilde{\mfD}_k)}+\|\bl{h}_k\|_{H^{\frac12}(\p\tilde{\mfD}_k)}\leq C_{\al,\mD}\|\bl{v}\|_{H^1(\mfD_k)}\left(1+\|\bl{v}\|_{H^1(\mfD_k)}\right)+C\Phi\chi_{|k|\leq1}.\\
\end{split}
\ee
From now on, Lemma \ref{LEM59}, with $m=0$, is valid for the system \eqref{cutEq1}. Combining \eqref{EEEST2} above and \eqref{ESTLEM28} in Lemma \ref{LEM59}, one arrives
\be\label{EESTT2}
\begin{split}
\|\bl{v}\|_{H^{2}(\mfD_k)}\leq& \|\bl{\mathfrak{u}}_k\|_{H^{2}(\tilde{\mfD}_k)}+\|\xi_k\|_{H^{3}(\tilde{\mfD}_k)}\\
\leq& C_{\al,\mD}\|\bl{v}\|_{H^{1}(\mfD_k)}\left(1+\|\bl{v}\|_{H^{1}(\mfD_k)}\right)+C\Phi\chi_{|k|\leq1}.
\end{split}
\ee
Now summing over $k\in\mathbb{Z}$, one proves \eqref{reg0} and \eqref{hesti} by a similar approach as we prove \eqref{ESTW13} before. The estimate of $\na\Pi$ in \eqref{hesti} follows directly from the equation \eqref{EQU1}$_1$ and \eqref{EESTT2} above. This completes the proof of Proposition \ref{PROP51}.

\qed

\begin{remark}\label{RMKK54}
In the proof of Proposition \ref{PROP51}, one notices that \eqref{EEEE1} and \eqref{EEEE11} can lead to the following uniform estimate of the pressure by summing over $k\in\mathbb{Z}$
\be\label{EEEE111}
\sum_{k\in\mathbb{Z}}\|\Pi-\overline{\Pi}_{\mfD_k}\|^2_{L^2(\mfD_k)}\leq C_{\al,\mD}\Phi^2<\i.
\ee
Moreover, since \eqref{EEEE111} is derived in the framework of the $H^1$-weak solution, it is valid for the case that $\p\mD$ is less regular.
\end{remark}

\qed

\subsection{Higher-order regularity and related estimates}
\q\ Following the route that we derive the solution belongs to $H^2$, now we are ready to derive higher-order regularity of $(\bl{u},p)$ via bootstrapping.

\begin{proposition}\label{threg}
Let $(\bl{u},p)$ be the weak solution to the problem \eqref{GL1}--\eqref{GL3}. Then
\[
(\bl{u}\,,p)\in C^{\infty}(\overline{\mD})\,.
\]
Meanwhile
\[
\bl{v}=\bl{u}-\bl{a},\q \Pi=p+\f{\Phi\int_{-\i}^{x_3}\eta(s)ds}{C_{P,R}}-\f{\Phi\int_{-\i}^{-x_3}\eta(s)ds}{C_{P,L}}
\]
satisfies
\be\label{hestik}
\|\bl{v}\|_{H^{m+2}(\mD)}+\|\nabla \Pi\|_{H^{m}(\mD)}\leq C_{m,\al,\mD}\Phi.
\ee
\end{proposition}
\pf
The proof follows from an induction argument. First the case of $m=0$ is already shown in Proposition \ref{PROP51}. Once the regularity estimate \eqref{hestik} is achieved with order $m\geq0$, one deduces that
\be\label{ESTt2}
\begin{split}
\left\|\nabla^{m+1}(\bl{v}\cdot\nabla{\bl{v}})\right\|_{L^2(\mD)}\leq C_{m,\mD}\left(\|\bl{v}\|_{L^\i(\mD)}\|\nabla^{m+2}\bl{v}\|_{L^2(\mD)}+\|\bl{v}\|_{W^{m+1,4}(\mD)}^2\right)\leq C_{m,\mD}\|\bl{v}\|_{H^{m+2}(\mD)}^2.
\end{split}
\ee
Therefore the Navier-Stokes system \eqref{NS}--\eqref{NBC}, which is equivalent to
\be\label{EQU31}
\left\{
\begin{array}{*{2}{ll}}
-\Dl \bl{v}+\na \Pi=\bl{g}, \q\na\cdot \bl{v}=0,\q\q&\text{in }\mD;\\[2mm]
2(\mathbb{S}\bl{v}\cdot\bl{n})_{\mathrm{tan}}+\al \bl{v}_{\mathrm{tan}}=0,\q \bl{v}\cdot\bl{n}=0, \q\q&\text{on }\p\mD,
\end{array}
\right.
\ee
where
\[
\bl{g}=-\left(\bl{v}\cdot\nabla \bl{v}+\bl{v}\cdot\nabla \bl{a}+\bl{a}\cdot\nabla \bl{v}+\bl{a}\cdot\nabla \bl{a}\right)+\Dl \bl{a}+\Phi\left(\f{\eta(x_3)}{C_{P,R}}+\f{\eta(-x_3)}{C_{P,L}}\right)\bl{e_3},
\]
enjoys
\[
\|\bl{g}\|_{H^{m+1}(\mD)}\leq C_{m,\al,\mD}\Phi<\i
\]
by direct calculations. Meanwhile, the problem \eqref{AUXHAR}
\[
\left\{
\begin{aligned}
&\Dl\xi_k=v_3\psi_k',\q\text{in }\tilde{\mfD}_k,\\
&\f{\p\xi_k}{\p\bl{n}}=0,\hskip .9cm\text{on }\p\tilde{\mfD}_k,\\
&\int_{\tilde{\mfD}_k}\xi_k(x)dx=0,
\end{aligned}
\right.
\]
now admits a unique solution in $H^{m+4}(\tilde{\mfD}_k)$ that satisfies
\be\label{ESTt3}
\|\xi_k\|_{H^{m+4}(\tilde{\mfD}_k)}\leq C_{m,\mD}\|\bl{v}\|_{H^{m+2}(\mfD_k)}.
\ee
Here the constant $C_{m,\mD}$ is independent with $k$, because every $\tilde{\mfD}_k$ is congruent with an element in $\left\{\tilde{\mfD}_{-1},\,\tilde{\mfD}_{0},\,\tilde{\mfD}_{1}\right\}$. Recalling the construction of \eqref{cutEq1}, we conclude $\bl{\mathfrak{u}}_k:=\psi_k\bl{v}-\nabla\xi_k$ satisfies
\[
\left\{
\begin{array}{*{2}{ll}}
-\Dl\bl{\mathfrak{u}}_k+\nabla\left(\psi_k(\Pi-\overline{\Pi}_{\mfD_k})\right)=\bl{g}_k,&\q\q\text{in }\tilde{\mfD}_k;\\[2mm]
\nabla\cdot\bl{\mathfrak{u}}_k=0,&\q\q\text{in }\tilde{\mfD}_k;\\[2mm]
2(\mathbb{S}\bl{\mathfrak{u}}_k\cdot\bl{n})_{\mathrm{tan}}+\al (\bl{\mathfrak{u}}_k)_{\mathrm{tan}}=\bl{h}_k,\q \bl{\mathfrak{u}}_k\cdot\bl{n}=0, &\q\q\text{on }\p\tilde{\mfD}_k,
\end{array}
\right.
\]
with
\[
\begin{split}
\bl{g}_k&:=\mathrm{div}\,\left(\psi_k\bl{F}+\nabla^2\xi_k\right)+\psi_k\bl{f}+\left(\Pi-\overline{\Pi}_{\mfD_k}\right)\psi_k'\bl{e_3}-2\psi_k'\p_{x_3}\bl{v}-\psi_k''v_3\bl{e_3}-\psi_k'\bl{F}\cdot\bl{e_3},\\
\bl{h}_k&:=-2\left(\left(\mathbb{S}\nabla\xi_k\right)\cdot\bl{n}\right)_{\mathrm{tan}}-\al\left(\na\xi_k\right)_{\mathrm{tan}}.
\end{split}
\]
By induction, together with \eqref{ESTt2} and \eqref{ESTt3}, one deduces that
\be\label{EEEEE}
\|\bl{g}_k\|_{H^{m+1}(\tilde{\mfD}_k)}+\|\bl{h}_k\|_{H^{m+\frac32}(\p\tilde{\mfD}_k)}\leq C_{m,\al,\mD}\|\bl{v}\|_{H^{m+2}(\mfD_k)}\left(1+\|\bl{v}\|_{H^{m+2}(\mfD_k)}\right)+C_m\Phi\chi_{|k|\leq1}
\ee
by the approach in the proof of Lemma \ref{LEMCOR}. Using the higher-order regularity for linear Stokes equations in Lemma \ref{LEM59}, together with estimates \eqref{ESTt3} and \eqref{EEEEE}, one proves the following $H^{m+3}$ estimate in $\tilde{\mfD}_k$
\[
\begin{split}
\|\bl{v}\|_{H^{m+3}(\mfD_k)}\leq& \|\bl{\mathfrak{u}}_k\|_{H^{m+3}(\tilde{\mfD}_k)}+\|\xi_k\|_{H^{m+4}(\tilde{\mfD}_k)}\\
\leq& C_{m,\al,\mD}\|\bl{v}\|_{H^{m+2}(\mfD_k)}\left(1+\|\bl{v}\|_{H^{m+2}(\mfD_k)}\right)+C_m\Phi\chi_{|k|\leq1}.
\end{split}
\]
Summing over $k\in\mathbb{Z}$, one concludes
\[
\|\bl{v}\|_{H^{m+3}(\mD)}\leq C_{\al,m,\mD}\Phi.
\]
Finally, the estimate of $\na\Pi$ in \eqref{hestik} follows directly from the equation \eqref{EQU31}$_1$ and the estimate above. This completes the proof of Proposition \ref{threg}, which indicates the validity of \eqref{HODEST} in Theorem \ref{THM16}.

\qed

\subsection{Exponential Decay of the weak solution}
\q\ In this subsection, we will show the ${H^1}$-norm exponential decay property of the solution. Our proof is carried out under the framework of the $H^1$-weak solution, which means, we only assume the solution satisfies the estimate in Theorem \ref{THM3.8}. However, with the help of the higher-order uniform estimates of the solution in Proposition \ref{threg}, the proof of the exponential decay property would be much simpler. Nevertheless, our proof in this subsection is also valid for stationary Navier-Stokes problem on domains which is less regular, say an infinite pipe only with a $C^{1,1}$ boundary.


\begin{proposition}\label{PROP5.1}
Let the conditions of Theorem \ref{Ext} be satisfied and $(\bl{v},\,\Pi)$ is given in \eqref{DEFVPI}. Then there exist positive constants $C$, $\sigma$, depending only on $\al$ and $\mD$, such that
\be\label{ASYP}
\begin{split}
\left\|\bl{u}-\bl{g}^L_{\Phi}\right\|_{H^1(\Sigma_L\times(-\i,-\zeta))}+\left\|\bl{u}-\bl{g}^R_{\Phi}\right\|_{H^1(\Sigma_R\times(\zeta,\i))}&\leq C\|\bl{v}\|_{H^1(\mD)}\exp(-\sigma\zeta),\\
\end{split}
\ee
for any $\zeta>Z+1$.
\end{proposition}
\pf
We only prove the estimate of term $\|\bl{u}-\bl{g}^R_{\Phi}\|_{H^1(\Sigma_R\times(\zeta,\infty))}$ since the rest term is essentially identical. In $\Sigma_R\times(Z,\i)$, the equation of $\bl{v}=\bl{u}-\bl{a}$ reads
\be\label{EW}
\bl{v}\cdot\nabla \bl{v}+\bl{a}\cdot\nabla\bl{v}+\bl{v}\cdot\na\bl{a}+\nabla\Pi-\Dl\bl{v}=0.
\ee
This is because
\[
\Dl \bl{a}+\left(\f{\eta(x_3)}{C_{P,R}}+\f{\eta(-x_3)}{C_{P,L}}\right)\Phi\bl{e_3}-\bl{a}\cdot\nabla\bl{a}=\left(\Dl g_\Phi^L+\f{\Phi}{C_{P,L}}\right)\bl{e_3}=0,\q\text{in}\q\Sigma_R\times(Z,\i).
\]
In the following proof, we will drop (upper or lower) indexes ``$R$" for convenience. For any $Z<\zeta\leq\zeta'<\zeta_1$, taking inner product with $\bl{v}$ on both sides of \eqref{EW} and integrating on $\Sigma\times(\zeta',\,\zeta_1)$, one has
\be\label{Maint000}
\un{\int_{\Sigma\times(\zeta',\,\zeta_1)}\bl{v}\cdot\Dl \bl{v}dx}_{LHS}=\un{\int_{\Sigma\times(\zeta',\,\zeta_1)}\big(\bl{v}\cdot\nabla \bl{v}+\bl{a}\cdot\nabla \bl{v}+\bl{v}\cdot\nabla \bl{a}+\nabla\Pi\big)\cdot\bl{v}dx}_{RHS}.
\ee
To handle the left hand side of \eqref{Maint000}, one first recalls the derivation of \eqref{Maint1} that
\be\label{SEC5M1}
\begin{split}
\int_{\Sigma\times(\zeta',\,\zeta_1)}\bl{v}\cdot\Dl\bl{v}dx=&-2\int_{\Sigma\times(\zeta',\,\zeta_1)}|\mathbb{S} \bl{v}|^2dx-\al\int_{\p\Sigma\times(\zeta',\,\zeta_1)}|\bl{v}_{tan}|^2dS\\
&-\sum_{i=1}^3\int_{\Sigma\times\{x_3=\zeta'\}}v_i(\p_{x_3}v_i+\p_{x_i}v_3)dx_h+\sum_{i=1}^3\int_{\Sigma\times\{x_3=\zeta_1\}}v_i(\p_{x_3}v_i+\p_{x_i}v_3)dx_h.
\end{split}
\ee
On the other hand, one can derive the following similarly as \eqref{Maint2}:
\be\label{SEC5M2}
\begin{split}
\int_{\Sigma\times(\zeta',\,\zeta_1)}\bl{v}\cdot\Dl \bl{v}dx\leq&-\int_{\Sigma\times(\zeta',\,\zeta_1)}|\nabla \bl{v}|^2dx+C_{\al,\mD}\int_{\p\Sigma\times(\zeta',\,\zeta_1)}|\bl{v}_{tan}|^2dS\\
&-\sum_{i=1}^3\int_{\Sigma\times\{x_3=\zeta'\}}v_i\p_{x_3}v_idx_h+\sum_{i=1}^3\int_{\Sigma\times\{x_3=\zeta_1\}}v_i\p_{x_3}v_idx_h.
\end{split}
\ee
Therefore, by calculating
\[
\eqref{SEC5M1}\times C_{\al,\mD}+\eqref{SEC5M2}\times\al,
\]
one deduces the left hand side of \eqref{Maint000} satisfies
\be\label{LH}
\begin{split}
LHS\leq&-\al\int_{\Sigma\times(\zeta',\,\zeta_1)}|\na \bl{v}|^2dx+C_{\al,\mD}\left(\sum_{i=1}^3\int_{\Sigma\times\{x_3=\zeta_1\}}v_i\p_{x_i}v_3dx_h-\sum_{i=1}^3\int_{\Sigma\times\{x_3=\zeta'\}}v_i\p_{x_i}v_3dx_h\right)\\
&+\left(\al+C_{\al,\mD}\right)\left(\sum_{i=1}^3\int_{\Sigma\times\{x_3=\zeta_1\}}v_i\p_{x_3}v_idx_h-\sum_{i=1}^3\int_{\Sigma\times\{x_3=\zeta'\}}v_i\p_{x_3}v_idx_h\right).
\end{split}
\ee
Using integration by parts for the right hand side of \eqref{Maint000}, one arrives
\be\label{RH}
\begin{split}
RHS=&\int_{\Sigma\times\{x_3=\zeta_1\}}\left(\f{1}{2}\left(v_3+g_\Phi\right)|\bl{v}|^2+v_3\Pi+g_\Phi(v_3)^2\right)dx_h\\
&-\int_{\Sigma\times\{x_3=\zeta'\}}\left(\f{1}{2}\left(v_3+g_\Phi\right)|\bl{v}|^2+v_3\Pi+g_\Phi(v_3)^2\right)dx_h\\
&-\int_{\Sigma\times(\zeta',\,\zeta_1)}\bl{v}\cdot\na \bl{v}\cdot \bl{a} dx.
\end{split}
\ee
Now we are ready to perform $\zeta_1\to\i$. To do this, one must be careful with the integrations on $\Sigma\times\{x_3=\zeta_1\}$ in both \eqref{LH} and \eqref{RH}. Recalling estimates of $(\bl{v},\Pi)$ in Theorem \ref{THM3.8} and Remark \ref{RMKK54}, one derives
\be\label{WESTTT}
\|\bl{v}\|^2_{H^1(\mD)}+\|\bl{v}\|^4_{L^4(\mD)}+\sum_{k\in\mathbb{Z}}\|\Pi-\overline{\Pi}_{\mfD_k}\|^2_{L^2(\mfD_k)}\leq C_{\al,\mD}\Phi^2<\i.
\ee
Choosing $M:=\f{C_{\al,\mD}\Phi^2}{Z}$, one concludes that for any $k>1$, there exists a slice $\Sigma\times\{x_3=\zeta_{1,k}\}$ which satisfies
\[
\Sigma\times\{x_3=\zeta_{1,k}\}\subset\mD\cap\left\{x\in\mathbb{R}^3:\,\left(\f{3k}{2}-\f{1}{2}\right)Z\leq x_3\leq\left(\f{3k}{2}+\f{1}{2}\right)Z\right\}\subset\mfD_k,
\]
and it holds that
\[
\int_{\Sigma\times\{x_3=\zeta_{1,k}\}}\left(|\na\bl{v}|^2+|\bl{v}|^4+|\Pi-\overline{\Pi}_{\mfD_k}|^2\right)dx_h\leq M.
\]
Otherwise, one has
\[
\|\bl{v}\|^2_{H^1(\mfD_k)}+\|\bl{v}\|^4_{L^4(\mfD_k)}+\|\Pi-\overline{\Pi}_{\mfD_k}\|^2_{L^2(\mfD_k)}>ZM=C_{\al,\mD}\Phi^2, \]
which creates a paradox to \eqref{WESTTT}. Choosing $k_0>0$ being sufficiently large such that the sequence $\{\zeta_{1,k}\}_{k=k_0}^\i\subset [\zeta',\i)$, clearly one has $\zeta_{1,k}\nearrow\i$ as $k\to\i$. Moreover, using the trace theorem of functions in the Sobolev space $H^1$, one has
\[
\int_{\Sigma}|\bl{v}(x_h,x_3)|^2dx_h\leq C\int_{z>x_3}\int_{\Sigma}(|\bl{v}|^2+|\nabla \bl{v}|^2)(x_h,z)dx_hdz\to 0,\q\text{as}\q x_3\to\i.
\]
Noting that $\int_{\Sigma\times\{x_3=\zeta_{1,k}\}}v_3dx_h=0$ for $k\geq k_0$, we deduce the following by the Poincar\'e inequality:
\[
\begin{split}
\left|\int_{\Sigma\times\{x_3=\zeta_{1,k}\}}v_3\Pi dx_h\right|&=\left|\int_{\Sigma\times\{x_3=\zeta_{1,k}\}}v_3\left(\Pi-\overline{\Pi}_{\mfD_k}\right) dx_h\right|\\
&\leq \left(\int_{\Sigma\times\{x_3=\zeta_{1,k}\}}|\bl{v}|^2dx_h\right)^{1/2}\left(\int_{\Sigma\times\{x_3=\zeta_{1,k}\}}|\Pi-\overline{\Pi}_{\mfD_k}|^2dx_h\right)^{1/2}\to 0,\q\text{as}\q k\to\i.
\end{split}
\]
Meanwhile, one finds
\[
\begin{split}
&\int_{\Sigma\times\{x_3=\zeta_{1,k}\}}|\bl{v}|\left(|\nabla\bl{v}|+|\bl{v}|^2\right)dx_h\\
&\leq \left(\int_{\Sigma\times\{x_3=\zeta_{1,k}\}}\left(|\na\bl{v}|^2+|\bl{v}|^4\right)dx_h\right)^{1/2}\left(\int_{\Sigma\times\{x_3=\zeta_{1,k}\}}|\bl{v}|^2dx_h\right)^{1/2}\to 0,\q\text{as}\q k\to\i;
\end{split}
\]
and
\[
\int_{\Sigma\times\{x_3=\zeta_{1,k}\}}|g_\Phi||\bl{v}|^2dx_h\leq\|g_\Phi\|_{L^\i(\mD_R)}\int_{\Sigma\times\{x_3=\zeta_{1,k}\}}|\bl{v}|^2dx_h\to0,\q\text{as}\q k\to\i.
\]
Choosing $\zeta_1=\zeta_{1,k}$ ($k\geq k_0$) in \eqref{LH} and \eqref{RH}, respectively, and performing $k\to\i$, one can deduce that
\[
\begin{split}
\al\int_{\Sigma\times(\zeta',\i)}|\na \bl{v}|^2dx\leq&\un{\int_{\Sigma\times(\zeta',\i)}\bl{v}\cdot\na \bl{v}\cdot \bl{a} dx}_{R_1}\\
&+C_{\al,\mD}\int_{\Sigma\times\{x_3=\zeta'\}}\Big(|\bl{v}|\left(|\bl{v}|^2+|g_\Phi||\bl{v}|+|\nabla\bl{v}|\right)+v_3\Pi\Big)dx_h.
\end{split}
\]
Using the Cauchy-Schwarz inequality, the Poincar\'e inequality in Lemma \ref{P2}, and the construction of profile vector $\bl{a}$, one derives
\[
R_1\leq\|\bl{a}\|_{L^\i(\mD)}\left(\int_{\Sigma\times(\zeta',\i)}|\na \bl{v}|^2dx\right)^{1/2}\left(\int_{\Sigma\times(\zeta',\i)}|\bl{v}|^2dx\right)^{1/2}\leq C_0\Phi\int_{\Sigma\times(\zeta',\i)}|\na \bl{v}|^2dx,
\]
which indicates the following estimate provided $\Phi$ is small enough such that $C_0\Phi<\al$:
\be\label{MEST}
\int_{\Sigma\times(\zeta',\i)}|\na \bl{v}|^2dx\leq C_{\al,\mD}\int_{\Sigma\times\{x_3=\zeta'\}}\Big(|\bl{v}|\left(|\bl{v}|^2+|g_\Phi||\bl{v}|+|\nabla\bl{v}|\right)+v_3\Pi\Big)dx_h.
\ee
Denoting
\be\label{DEFG}
\mathcal{G}(\zeta'):=\int_{\Sigma\times(\zeta',\i)}|\na \bl{v}|^2dx,
\ee
and integrating \eqref{MEST} with $\zeta'$ on $(\zeta,\i)$, one arrives
\be\label{EEP0}
\int_\zeta^\i\mathcal{G}(\zeta')d\zeta'\leq C_{\al,\mD}\left(\int_{\Sigma\times(\zeta,\i)}\Big(|\bl{v}|\left(|\bl{v}|^2+|g_\Phi||\bl{v}|+|\nabla\bl{v}|\right)\Big)dx+\left|\int_{\Sigma\times(\zeta,\i)}v_3\Pi dx\right|\right).
\ee
Applying the Poincar\'e inequality in Lemma \ref{P2}, one deduces
\be\label{EEP1}
\int_{\Sigma\times(\zeta,\i)}|\bl{v}|\left(|\bl{v}|^2+|g_\Phi||\bl{v}|+|\nabla\bl{v}|\right)dx\leq C_{\al,\mD}\int_{\Sigma\times(\zeta,\i)}|\nabla\bl{v}|^2dx.
\ee
Moreover, using a similar approach as in the proof of Proposition \ref{P2.4}, one notices that
\be\label{EEP2}
\begin{split}
\left|\int_{\Sigma\times(\zeta,\i)}v_3\Pi dx\right|&\leq\sum_{m=1}^\i\left|\int_{\O_{\zeta+m}^+}v_3\Pi dx\right|\\
&\leq C\sum_{m=1}^\i\left(\|g_\Phi\|_{L^\infty(\O^{+}_{\zeta+m})}\|\nabla \bl{v}\|^2_{L^2(\O^{+}_{\zeta+m})}+\|\nabla \bl{v}\|_{L^2(\O^{+}_{\zeta+m})}^2+\|\nabla \bl{v}\|_{L^2(\O^{+}_{\zeta+m})}^3\right)\\
&\leq C\int_{\Sigma\times(\zeta,\i)}|\nabla\bl{v}|^2dx.
\end{split}
\ee
Substituting \eqref{EEP1} and \eqref{EEP2} in \eqref{EEP0}, one arrives at
\[
\int_\zeta^\i\mathcal{G}(\zeta')d\zeta'\leq C_{\al,\mD}\mathcal{G}(\zeta),\q\text{for any}\q\zeta>Z.
\]
This implies
\[
\mathcal{N}(\zeta):=\int_\zeta^\i\mathcal{G}(\zeta')d\zeta'
\]
is well-defined for all $\zeta>Z$, and
\be\label{EEP3}
\mathcal{N}(\zeta)\leq-C_{\al,\mD}\mathcal{N}'(\zeta),\q\text{for any}\q\zeta>Z.
\ee
Multiplying the factor $e^{C^{-1}_{\al,\mD}\zeta}$ on both sides of \eqref{EEP3} and integrating on $[Z,\zeta]$, one deduces
\[
\mathcal{N}(\zeta)\leq C_{\al,\mD}\exp\left(-C^{-1}_{\al,\mD}\zeta\right),\q\text{for any}\q\zeta>Z.
\]
According to the definition \eqref{DEFG}, one has $\mathcal{G}$ is both non-negative and non-increasing. Thus
\[
\mathcal{G}(\zeta)\leq\int_{\zeta-1}^\zeta\mathcal{G}(\zeta')d\zeta'\leq\mathcal{N}(\zeta-1)\leq C\exp\left(-C^{-1}_{\al,\mD}\zeta\right),\q\text{for any}\q\zeta>Z+1.
\]
This completes the proof of the \eqref{ASYP} by choosing $\sigma=C^{-1}_{\al,\mD}$.

\qed

\subsection{On the exponential decay for higher-order derivatives}

\q\ In this section, we focus on the higher-order asymptotic behavior of the aforementioned unique smooth solution to the problem. One will see that the solution $\bl{u}$ converges to the Poiseuille flow at the far field with an exponential speed. Based on the $H^1$ decay property in Proposition \ref{PROP5.1}, we finish the proof of the estimate \eqref{pointdecay} in Theorem \ref{THM16}.

Using Sobolev imbedding, we first need to show the following decay of the solution in $H^m$ norms, with $m\geq 2$:
\[
\|\bl{v}\|_{H^m(\Sigma_L\times(-\i,-\zeta))}+\|\bl{v}\|_{H^m(\Sigma_R\times(\zeta,\i))}\leq C_{m,\al,\mD}\left(\|\bl{v}\|_{H^1(\Sigma_L\times(-\i\,,\,-\zeta+2Z))}+\|\bl{v}\|_{H^1(\Sigma_R\times(\zeta-2Z\,,\,\i))}\right),\q\text{for all }\zeta>3Z.
\]
This is derived by using the method in the proof of Proposition \ref{threg}, but summing over $k\in\mathbb{Z}$ such that
\[
\mathrm{supp}\,\psi_k\cap\big((-\i,-\zeta)\cup(\zeta,\i)\big)\neq\varnothing.
\]
Then, the proof is completed by the $H^1$ decay estimate \eqref{ASYP}.

\section*{Data availability statement}
\addcontentsline{toc}{section}{Data availability statement}
\q\ Data sharing is not applicable to this article as no datasets were generated or analysed during the current study.

\section*{Conflict of interest statement}
\addcontentsline{toc}{section}{Conflict of interest statement}
\q\ The authors declare that they have no conflict of interest.

\section*{Acknowledgments}
\addcontentsline{toc}{section}{Acknowledgments}
\q\ Z. Li is supported by Natural Science Foundation of Jiangsu Province (No. BK20200803) and National Natural Science Foundation of China (No. 12001285). X. Pan is supported by National Natural Science Foundation of China (No. 11801268, 12031006). J. Yang is supported by National Natural Science Foundation of China (No. 12001429).

\medskip
\medskip

 {\footnotesize

{\sc Z. Li: School of Mathematics and Statistics, Nanjing University of Information Science and Technology, Nanjing 210044, China}

  {\it E-mail address:}  zijinli@nuist.edu.cn

\medskip

 {\sc X. Pan: College of Mathematics and Key Laboratory of MIIT, Nanjing University of Aeronautics and Astronautics, Nanjing 211106, China}

  {\it E-mail address:}  xinghong\_87@nuaa.edu.cn

\medskip

 {\sc J. Yang: School of Mathematics and Statistics, Northwestern Polytechnical University, Xi'an 710129, China}

  {\it E-mail address:} yjqmath@nwpu.edu.cn
}
\end{document}